\documentclass[]{elsart}
\newcommand{\diag}{{\rm diag}} 
\usepackage{times}
\usepackage{amssymb}
\usepackage{amsmath,fancyhdr}
\usepackage{graphicx}
\graphicspath{../}
\usepackage{url}
\usepackage{color}
\newcommand{\wigner}[6] 
{\left(\barray{ccc} #1 & #2 & #3\\ #4 & #5 & #6 \earray\right)}
\newcommand{\fracd}[2]{\frac{\displaystyle{#1}}{\displaystyle{#2}}}
\newcommand{\nnr}{\nonumber}
\newcommand{\barray}{\begin{array}}
\newcommand{\earray}{\end{array}}
\newcommand{\ber}{\begin{eqnarray}}
\newcommand{\eer}{\end{eqnarray}}
\newcommand{\T}{^{\sf{T}}}
\newcommand{\hsom}{\hspace{0.1em}}
\newcommand{\hsomm}{\hspace{-0.1em}}
\newcommand{\where}{\qquad\mbox{where}\qquad}
\newcommand{\also}{\qquad\mbox{and}\qquad}

\newcommand{\Xlmth}{X_{lm}(\theta)}
\newcommand{\bu}{{\bf u}}
\newcommand{\su}{\mathsf{u}}
\newcommand{\bv}{{\bf v}}

\newcommand{\but}{{\bf u}^t}
\newcommand{\bur}{{\bf u}^r}
\newcommand{\bgt}{{\bf g}^t}
\newcommand{\bgr}{{\bf g}^r}
\newcommand{\sg}{\mathsf{g}}
\newcommand{\sh}{\mathsf{h}}
\newcommand{\bg}{{\bf g}}
\newcommand{\bh}{{\bf h}}

\newcommand{\sB}{\mathsf{B}} %

\newcommand{\sP}{\mathsf{P}} %
\newcommand{\sQ}{\mathsf{Q}} %
\newcommand{\sC}{\mathsf{C}} %
\newcommand{\sD}{\mathsf{D}} %
\newcommand{\sK}{\mathsf{K}} %
\newcommand{\so}{\mathsf{0}}
\newcommand{\sPm}{\sP_m}
\newcommand{\sBm}{\sB_m}
\newcommand{\sDm}{\sD_m}
\newcommand{\sQm}{\sQ_m}
\newcommand{\sBmm}{\sB_{-m}}
\newcommand{\sDmm}{\sD_{-m}}
\newcommand{\sPmm}{\sP_{-m}}
\newcommand{\Ulm}{U_{lm}}
\newcommand{\Ulmp}{U_{l'm'}}
\newcommand{\Vlm}{V_{lm}}
\newcommand{\Vlmp}{V_{l'm'}}
\newcommand{\Wlm}{W_{lm}}
\newcommand{\Wlmp}{W_{l'm'}}
\newcommand{\Xlm}{X_{lm}}
\newcommand{\Xlamth}{X_{l \lvert m \rvert}(\theta)}
\newcommand{\Xlpm}{X_{l'm}}

\newcommand{\Xlpmpth}{X_{l'm'}(\theta)} %

\newcommand{\Xloth}{X_{l0}(\theta)} %

\newcommand{\Ylpmp}{Y_{l'm'}}
\newcommand{\Ylm}{Y_{lm}} %
\newcommand{\deltamm}{\delta_{mm'}}
\newcommand{\deltall}{\delta_{ll'}}
\newcommand{\deltaab}{\delta_{\alpha\beta}}
\newcommand{\renormprod}{\sqrt{l(l+1)l'(l'+1)}} %
\newcommand{\renormsing}{\sqrt{l(l+1)}} %
\newcommand{\Plm}{{\bf P}_{lm}}
\newcommand{\Coo}{{\bf C}_{00}} %
\newcommand{\Poo}{{\bf P}_{00}} %
\newcommand{\Boo}{{\bf B}_{00}} %
\newcommand{\Plmp}{{\bf P}_{l'm'}}
\newcommand{\Blm}{{\bf B}_{lm}}
\newcommand{\Blmp}{{\bf B}_{l'm'}}
\newcommand{\Clm}{{\bf C}_{lm}}
\newcommand{\Clmp}{{\bf C}_{l'm'}}
\newcommand{\Blmlmp}{B_{lm,l'm'}} %
\newcommand{\Clmlmp}{C_{lm,l'm'}} %
\newcommand{\Plmlmp}{P_{lm,l'm'}} %
\newcommand{\Plmlpm}{P_{lm,l'm}} %
\newcommand{\Blmlpm}{B_{lm,l'm}} %
\newcommand{\Dlmlpm}{D_{lm,l'\,-m}} %
\newcommand{\Pllp}{P_{ll'}} %
\newcommand{\Blmlm}{B_{lm,lm}} %
\newcommand{\Clmlm}{C_{lm,lm}} %
\newcommand{\Plmlm}{P_{lm,lm}} %
\newcommand{\Bllp}{B_{ll'}} 
\newcommand{\Dllp}{D_{ll'}} %
\newcommand{\Dlmlmp}{D_{lm,l'm'}} %
\newcommand{\bnabla}{\mbox{\boldmath$\nabla$}}
\newcommand{\Pll}{P_{ll}} %
\newcommand{\Bll}{B_{ll}} %

\newcommand{\gauntsum}{\sum_{n=\lvert l-l'\rvert}^{l+l'}} 
\newcommand{\stdsumL}{\sum_{lm}^L}
\newcommand{\stdsum}{\sum_{lm}^\infty}
\newcommand{\stdsumLp}{\sum_{l'm'}^L}

\newcommand{\aplm}{a^+_{lm}}
\newcommand{\amlm}{a^-_{lm}}

\newcommand{\apmlm}{a^\pm_{lm}}
\newcommand{\aplpm}{a^+_{l'm}}
\newcommand{\amlpm}{a^-_{l'm}}
\newcommand{\bplm}{b^+_{lm}}
\newcommand{\bmlm}{b^-_{lm}}

\newcommand{\bpmlm}{b^\pm_{lm}}
\newcommand{\bplpm}{b^+_{l'm}}
\newcommand{\bmlpm}{b^-_{l'm}}
\newcommand{\intOmega}{\int_{\Omega}}
\newcommand{\intoth}{\int_0^\Theta}
\newcommand{\intR}{\int_{R}}
\newcommand{\dOmega}{\,d\Omega}
\newcommand{\sindth}{\sin \theta \, d\theta}
\newcommand{\dXlm}{X'_{lm}}

\newcommand{\dXlpm}{X'_{l'm}}
\newcommand{\divsin}{(\sin \theta)^{-1}}
\newcommand{\divsinsq}{(\sin \theta)^{-2}}

\newcommand{\IlmTh}{I_{lm}(\Theta)}
\newcommand{\rvec}{{\bf \hat{r}}}
\newcommand{\thvec}{\boldsymbol{\hat \theta}}
\newcommand{\phvec}{\boldsymbol{\hat \phi}}
\newcommand{\ovec}{{\bf 0}}
\newcommand{\di}{\boldsymbol{\delta}}
\newcommand{\Di}{{\bf K}}

\newcommand{\dirrp}{\di(\rvec,\rvec')}

\newcommand{\Dirrp}{\Di(\rvec,\rvec')}

\newcommand{\bandlim}{\mathcal{S}_L}
\newcommand{\spacelim}{\mathcal{S}_R}
\newcommand{\sT}{\mathsf{T}}
\newcommand{\iT}{T}
\newcommand{\bga}{{\bf g}_{\alpha}}

\newcommand{\sga}{\mathsf{g}_{\alpha}}
\newcommand{\bgb}{{\bf g}_{\beta}}

\newcommand{\sgb}{\mathsf{g}_{\beta}}
\newcommand{\bha}{{\bf h}_{\alpha}}

\newcommand{\bhb}{{\bf h}_{\beta}}

\newcommand{\dphi}{\partial_\phi}
\newcommand{\dth}{\partial_\theta}
\newcommand{\dr}{\partial_r}
\newcommand{\sumalpha}{\sum_{\alpha=1}^{3(L+1)^2-2}}
\newcommand{\sumalphaJ}{\sum_{\alpha=1}^{J}}
\newcommand{\maxm}{\max(m,1)}

\newcommand{\bmp}{\begin{minipage}} 
\newcommand{\emp}{\end{minipage}} 
\newcommand{\edc}{\end{document}} 
\newcommand{\figend}{\\[1em]\hrule}
\DeclareMathOperator{\tr}{tr}

\begin{document}
\begin{frontmatter}
\title{Spatiospectral concentration of vector fields on a sphere}
\vspace{-0.75cm}
\author[PU]{Alain Plattner\corauthref{cor1}} and
\author[PU]{Frederik J.~Simons}
\corauth[cor1]{Corresponding author. E-mail: plattner@alumni.ethz.ch}
\address[PU]{Department of Geosciences, Princeton University, Princeton
 NJ 08544, U.S.A.}

\vspace{-2cm}
\received{\today}
\revised{\today}
\accepted{\today}

\lhead[]{}
\chead[]{\scriptsize\it{A.~Plattner and F.~J.~Simons /
   Appl.~Comput.~Harmon.~Anal. XX (2012) XXX-XXX}} 
\rhead[]{}

\vspace{2cm}

\begin{abstract}
  \hspace{0.4cm} We construct spherical vector bases that are
  bandlimited and spatially concentrated, or alternatively,
  spacelimited and spectrally concentrated, suitable for the analysis
  and representation of real-valued vector fields on the surface of
  the unit sphere, as arises in the natural and biomedical sciences,
  and engineering. Building on the original approach of Slepian,
  Landau, and Pollak we concentrate the energy of our function basis
  into arbitrarily shaped regions of interest on the sphere and within
  a certain bandlimit in the vector spherical-harmonic domain. As with
  the concentration problem for scalar functions on the sphere, which
  has been treated in detail elsewhere, the vector basis can be
  constructed by solving a finite-dimensional algebraic eigenvalue
  problem. The eigenvalue problem decouples into separate problems for
  the radial, and tangential components. For regions with advanced
  symmetry such as latitudinal polar caps, the spectral concentration
  kernel matrix is very easily calculated and block-diagonal, which
  lends itself to efficient diagonalization.  The number of
  spatiospectrally well-concentrated vector fields is well estimated
  by a Shannon number that only depends on the area of the target
  region and the maximal spherical harmonic degree or bandwidth. The
  Slepian spherical vector basis is doubly orthogonal, both over the
  entire sphere and over the geographic target regions. Like its
  scalar counterparts it should be a powerful tool in the inversion,
  approximation and extension of bandlimited fields on the sphere:
  vector fields such as gravity and magnetism in the earth and
  planetary sciences, or electromagnetic fields in optics, antenna
  theory and medical imaging.

\end{abstract}

\begin{keyword}
bandlimited function, concentration problem, eigenvalue
 problem, spectral analysis, vector spherical harmonics
\end{keyword}

\end{frontmatter}

 \section{Introduction}

 Since it is impossible to simultaneously bandlimit and spacelimit a
 function to a chosen region of interest, we need to resort to
 functions that are band\-limited but optimally concentrated, with
 respect to their spatial energy, inside a target region.  Slepian,
 Landau, and Pollak presented the solution for the problem of
 optimally concentrating a signal in time and frequency in their
 seminal papers
 \cite{Landau+61,Landau+62,Slepian83,Slepian+61}.  Their construction
 leads to a family of orthogonal tapers or data windows that have been
 widely applied as windows to regularize the quadratic inverse problem
 of power spectral estimation from time-series observations of finite
 extent. The ``Slepian functions'', as we shall be calling them, are
 furthermore of great utility as a basis for function representation,
 approximation, interpolation and extension, and to solve stochastic
 linear inverse problems in a wide range of disciplines.  Several
 authors have studied the time-scale and time-frequency concentration
 for more general settings
 (see~\cite{Simons2010,Simons+2006a,Simons+2011a} and references
 therein for a review).  More specifically, spherical scalar Slepian
 functions, spatially concentrated while bandlimited or spectrally
 concentrated while spacelimited, have been applied in physical,
 computational, and biomedical fields such as
 geodesy~\cite{Albertella2001,Simons+2006b,Han+2008b,Slobbe+2012} and
 gravimetry
  \cite{Han2008c,Han2008a,Han2008d,Longuevergne+2010,Harig+2012},
 geomagnetism~\cite{Schott+2011,Simons+2009b} and
 geodynamics~\cite{Harig2010},
 planetary~\cite{Evans+2010,Goossens+2012,Han2009,Wieczorek2008} and
 biomedical science~\cite{Maniar2005,Mitra2006},
 cosmology~\cite{Dahlen2008,Das+2009}, and computer
 science~\cite{Lessig+2010,Yeo+2008}, while continuing to be of
 interest in 
information and communication theory~\cite{SenGupta+2012}, signal
processing~\cite{Khalid+2012,Wei+2011}, and
mathematics~\cite{Marinucci+2010,Michel2011}.  

To date only a few attempts have been made to bring the advantages of
spherical Slepian functions into the realm of spherical vector
fields. The first successful construction of spatially concentrated
bandlimited tangential spherical vector fields was reported for
applications in
magnetoencephalography~\cite{Maniar2005,Mitra2006,Mitra+2005}. In
geodesy, 
Eshagh~\cite{Eshagh2009a} has developed methods to explicitly evaluate
the product integrals arising in the concentration problem whose solutions
are the vectorial Slepian functions. In this paper we present a
complete extension of Slepian's spatiospectral concentration problem
to vector fields on the sphere, and give suggestions and examples as
to their usage for problems of a geomagnetic nature
(e.g.~\cite{Stockmann+2009,Gubbins+2011b}).  
\begin{figure}[bt]
\begin{centering}
 \includegraphics[width=0.65\textwidth]{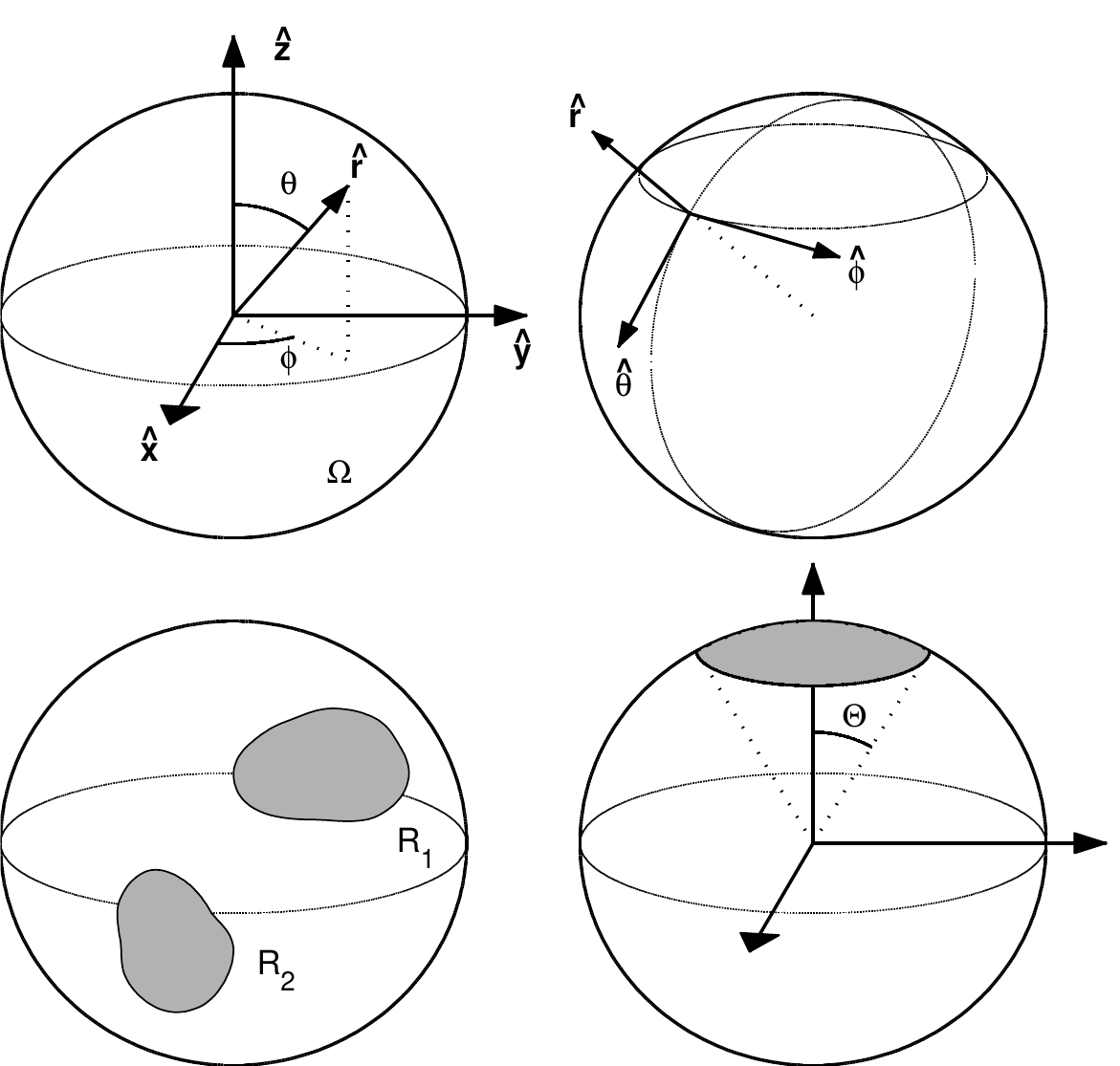}\figend
 \caption{\label{diagram} Sketch illustrating the geometry of the
   vector spherical concentration problem.  Lower right shows an
   axisymmetric polar cap of colatitudinal radius~$\Theta$ as treated
   in Section~\ref{polar cap section} 
   The area of the region of concentration, $R=R_1 \cup R_2 \cup
   \ldots$, is denoted by~$A$ in the text.}
\end{centering}
\end{figure}
The family of optimally concentrated spherical vectorial multitapers
that we will construct in the following should be useful in many
scientific applications.  In particular in geomagnetism, one of the
objectives of the \textsc{swarm} mission~\cite{Friis2006} is to model
the lithospheric magnetic field with maximal resolution and accuracy,
even in the presence of contaminating signals from secondary sources.
In addition, and more generally, lithospheric-field data analysis will
have to successfully merge information from the global to the regional
scale. In the past decade or so, a variety of global-to-regional
modeling techniques have come of age including harmonic
splines~\cite{Shure1982,Parker1982,Langel1987}, stitching together
local
models~\cite{Haines1985,DeSantis1991,Korte2003,Thebault2006,Thebault2006a},
and wavelets~\cite{Holschneider2003,Chambodut2005,Mayer2006}.  Due to
their optimal combination of spatial locality and spectral
bandlimitation the multitapers constructed in this paper should be
well suited to combine global and local data while respecting the
bandlimitation.

\section{Preliminaries}

Figure~\ref{diagram} shows the geometry of the unit sphere~$\Omega =
\{\rvec : \lVert \rvec \rVert=1\}$ and its tangential vectors. The
colatitude of spherical points $\rvec$ is denoted by $0\leq \theta\leq
\pi$ and the longitude by $0\leq \phi<2\pi$; we denote the unit vector
pointing outwards in the radial direction by~$\rvec$, and the unit
vectors in the tangential directions towards the south pole and
towards the east will be denoted by~$\thvec$ and~$\phvec$,
respectively. The symbol~$R$ will be used to denote a region of the
unit sphere~$\Omega$, of area~$A=\int_R\dOmega$, within which the
bandlimited vector field shall be concentrated.  The region can be a
combination of disjoint subregions, $R=R_1\cup R_2 \cup \cdots$, and
the boundaries of the $R_i$ can be irregularly shaped, as depicted. We
will denote the region complementary to~$R$ by~$\Omega\setminus R$.

\subsection{Real Scalar Spherical Harmonics}

Restricting our attention to real-valued vector fields, we use real
vector spherical harmonics, which are constructed from their scalar
counterparts. Each scalar spherical harmonic~$\Ylm$ has a degree
$0\leq l$ and, for each degree, an order $-l\leq m\leq l$. Our
spherical harmonics are unit-normalized in the sense~\cite{Dahlen1998}
\begin{align}
\Ylm(\theta,\phi)&=
\begin{cases}  
\sqrt{2}\Xlamth \cos m\phi &\text{if } -l\leq m<0,\\
\Xloth&\text{if } m=0,\\
\sqrt{2}X_{l m}(\theta)\sin m\phi &\text{if } 0<m\leq l,
\end{cases}\label{Y definition}\\
\Xlmth&=(-1)^m\left(\frac{2l+1}{4\pi}\right)^{1/2}
\left[\frac{(l-m)!}{(l+m)!}\right]^{1/2} P_{lm}(\cos \theta),\label{X
definition}\\ 
P_{lm}(\mu)&=\frac{1}{2^ll!}(1-\mu^2)^{m/2}\left(\frac{d}{d\mu}\right)^{l+m}(\mu^2-1)^l.\label{P definition}
\end{align}
The asymptotic wavenumber~\cite{Jeans1923} associated with a harmonic
degree~$l$ is $\sqrt{l(l+1)}$. The function~$P_{lm}(\mu)$ in (\ref{P
  definition}) is called the associated Legendre function of integer
degree~$l$ and order~$m$. The spherical harmonics $\Ylm(\rvec)$ are
eigenfunctions of the Laplace-Beltrami operator, $\nabla_1^2
\Ylm=-l(l+1)\Ylm$, where $\nabla_1^2=\partial_\theta^2 +
\cot{\theta}\,\partial_\theta + (\sin{\theta})^{-2}\partial_\phi^2$.
We choose the constants in (\ref{Y definition})--(\ref{P definition})
to guarantee orthonormality:
\begin{equation}
\intOmega \Ylm \Ylpmp \dOmega = \deltall\deltamm.
\end{equation}
The product of two normalized Legendre
functions~\cite{Dahlen1998,Edmonds1996} is the linear
combination
\begin{align}\label{Gaunt}
\nnr \Xlmth \Xlpmpth
=(-1)^{m+m'}\gauntsum  \sqrt{\frac{(2n+1)(2l+1)(2l'+1)} 
{4\pi}}\\
\times
\wigner{l}{n}{l'}{0}{0}{0}
\wigner{l}{n}{l'}{m}{-m-m'}{m'}
 X_{n\,m+m'}(\theta).
\end{align}
The index arrays in (\ref{Gaunt}) are Wigner~3-$j$~symbols
\cite{Edmonds1996,Messiah2000}. We will use the following two recursion
relations~\cite{Eshagh2009a,Ilk1983} for the derivatives of
the~$\Xlmth$ and their divisions by $\sin\theta$. Define 
\begin{equation}\label{remind}
\dXlm= \frac{dX_{lm}}{d\theta}.
\end{equation}
Then, from the work  by Ilk~\cite{Ilk1983} follows that, for $0\leq l$
and $0\leq m \leq l$, 
\begin{equation}\label{Ilk 1}
\dXlm=\amlm X_{l\,m-1} + \aplm X_{l\,m+1}, 
\end{equation}
where
\begin{equation}\label{ilk coefficients a}
  \apmlm = \pm\frac{\sqrt{(l\mp m)(l\pm m+1)}}{2}.
\end{equation}
and, for $0\leq l$ and $1\leq m \leq l$,
\begin{equation}\label{Ilk 2}
  \divsin m\Xlm =\bmlm X_{l-1\,m-1}+\bplm X_{l-1\,m+1},
\end{equation}
where
\begin{equation}\label{ilk coefficients b}
\bpmlm=-\sqrt{\frac{2l+1}{2l-1}}\frac{\sqrt{(l\mp m)(l\mp m-1)}}{2}.
\end{equation}
In both~(\ref{ilk coefficients a}) and~(\ref{ilk coefficients b}), for
$l<0$ or $m<0$, we have~$\Xlm=0$. 

Finally, as shown by Paul~\cite{Paul1978}, integrals of the type
\begin{equation}\label{bladipaul}
 \IlmTh=\intoth \Xlmth \sindth
\end{equation}
can be exactly evaluated recursively. When $l\geq 2$ and $0\leq m<l$,
we have
\begin{align}
 \nnr\IlmTh=&
 \frac{l-2}{l+1}\sqrt{\frac{(2l+1)[(l-1)^2-m^2]}{(2l-3)(l^2-m^2)}}
 I_{l-2\, m}(\Theta)\\& +
 \frac{1}{l+1}\sqrt{\frac{4l^2-1}{l^2-m^2}}\sin^2\hsomm\Theta
 X_{l-1\,m}(\Theta)\label{Paul lm},   
\end{align}
and, for $l\geq 2$ and $m=l$, the formula
\begin{align}
 I_{ll}(\Theta)=&\frac{1}{l+1}\sqrt{\frac{2l+1}{4l^2-4l}}
 \left[l\sqrt{2l-1} I_{l-2\, l-2}(\Theta) -(\sin\Theta)^2
   X_{l-1\,l-2}(\Theta)\right]\label{Paul ll}  
.
\end{align}
The recursions (\ref{Paul lm})--(\ref{Paul ll}) are to be started from 
\begin{align}
 X_{00}(\Theta)&= \frac{1}{2\sqrt{\pi}},
 \qquad &I_{00}(\Theta)&=\frac{1}{2\sqrt{\pi}}(1-\cos\Theta),\\
 X_{10}(\Theta)&= \frac{1}{2}\sqrt{\frac{3}{\pi}}\cos \Theta, \qquad
 &I_{10}(\Theta)&=\frac{1}{4}\sqrt{\frac{3}{\pi}}(\sin\Theta)^2,\\ 
 X_{11}(\Theta)&=-\frac{1}{2}\sqrt{\frac{3}{2\pi}}\sin\Theta,  \qquad
 &I_{11}(\Theta)&=-\frac{1}{4}\sqrt{\frac{3}{2\pi}}(\Theta-\sin2\Theta),
\end{align}
hence enabling the exact evaluation of all integrals $I_{lm}(\Theta)$ for $l\geq 0$
and $0\leq m \leq l$.

\subsection{Real Valued Vector Spherical Harmonics}

The canonical three-dimensional gradient operator $\bnabla = \hat{\bf
  x}\hsom\partial_{x} + \hat{\bf y}\hsom\partial_{y} + \hat{\bf z}\hsom\partial_{z}$ 
can be expressed as~\cite{Dahlen1998}
\begin{equation}
 \bnabla= \rvec\hsom\dr + r^{-1}\bnabla_1,
\where
 \bnabla_1 = \thvec\hsom\dth + \phvec \hsom\divsin \dphi.
\end{equation}
For any differentiable function $H(\rvec)$ on the unit sphere, the
vector field $\rvec H(\rvec)$ is purely radial, the vector fields
$\bnabla_1 H(\rvec)$ and $\rvec\times \bnabla_1 H(\rvec)$ are purely
tangential, and all three are mutually orthogonal. We can thus
construct vector spherical harmonics from gradients of scalar
spherical harmonics by defining, for $l>0$ and $-m\leq l \leq m$,
\begin{align}
\Plm& = \rvec \Ylm,\label{Plm}\\
\Blm& = \frac{ \bnabla_1 \Ylm}{\sqrt{l(l+1)}} = \frac{[
  \thvec\hsom\dth + \phvec\hsom\divsin \dphi]
  \Ylm}{\renormsing},\label{Blm}\\ 
\Clm& = \frac{-\rvec\times \bnabla_1
  \Ylm}{\sqrt{l(l+1)}}=\frac{[\thvec\hsom\divsin \dphi-\phvec\hsom\dth]
  \Ylm}{\renormsing}, \label{Clm}
\end{align}
together with the purely radial $\Poo=(4\pi)^{-1/2}\hsom\rvec$ and the vanishing
$\Boo=\Coo=\ovec$. The orthonormality of $\rvec$, $\thvec$, and $\phvec$
immediately leads to \begin{equation}\label{strong orthogonality}
 \Plm\cdot\Blmp= \Plm\cdot\Clmp=\ovec,
\end{equation}
and they are furthermore orthonormal in the sense
\begin{align}\label{orthonormality}
 \intOmega \Plm\cdot\Plmp\dOmega&=\intOmega
 \Blm\cdot\Blmp\dOmega=\intOmega
 \Clm\cdot\Clmp\dOmega=\deltall\deltamm,\\ \label{orthonormality 2} 
 \intOmega \Plm\cdot\Blmp\dOmega&=\intOmega
 \Plm\cdot\Clmp\dOmega=\intOmega \Blm\cdot\Clmp\dOmega=0 
.\end{align}
The vector spherical-harmonic addition theorem~\cite{Freeden2009}
comprises the identities
\begin{eqnarray}
\label{vecadd1}\label{vecadd2}
\lefteqn{\sum_{m=-l}^l \Plm(\rvec)\cdot \Plm(\rvec)=
\left(\frac{2l+1}{4\pi}\right)}\\
&&
\hspace{3.5em}=\label{vecadd3}
\sum_{m=-l}^l \Blm(\rvec)\cdot \Blm(\rvec)=
\sum_{m=-l}^l \Clm(\rvec)\cdot \Clm(\rvec).
\end{eqnarray}

\subsection{Real-Valued Vector Fields on the Unit Sphere}

The expansion of a real-valued square-integrable vector
field~$\bu$ on the unit sphere~$\Omega$ can be written as 
\begin{equation}\label{general representation}
 {\bf u}= \stdsum \Ulm\Plm+ \Vlm\Blm + \Wlm\Clm,
\end{equation}
where the expansion coefficients are obtained via 
\begin{equation}
 \Ulm=\intOmega \Plm \cdot \bu\dOmega,\label{UlmVlmWlm}\quad
 \Vlm=\intOmega \Blm \cdot \bu\dOmega,\quad\mbox{and}\,\,
 \Wlm=\intOmega \Clm \cdot \bu\dOmega,
\end{equation}
using the shorthand notation $\stdsumL:=\sum_{l=0}^L\sum_{m=-l}^l$
when $\Plm$ or $\Ulm$ are involved and
$\stdsumL:=\sum_{l=1}^L\sum_{m=-l}^l$, for $\Blm$, $\Clm$, $\Vlm$ or
$\Wlm$.  A sans serif~$\su$ will be used to denote the ordered column
vector of vector spherical-harmonic coefficients, namely $\su=(\ldots,
\Ulm,\ldots,\Vlm,\ldots,\Wlm,\ldots)\T$.  We will denote the norms of
a spatial-domain vector field $\bu(\rvec)$ and its spectral-domain
equivalent $\su$ by
\begin{equation}\label{norms}
 \lVert \bu\rVert^2_\Omega=\intOmega \bu\cdot\bu\dOmega, \qquad
 \lVert \su\rVert^2_\infty= \sum_{lm}^\infty \Ulm^2 + \Vlm ^2+
 \Wlm^2. 
\end{equation}
Hence Parseval's relation can be written in the form
$\lVert \bu\rVert^2_\Omega = \lVert \su\rVert^2_\infty$. Any
square-integrable vector field $\bu$ on the sphere can be decomposed
into a radial component, $\bur$, and a tangential component, $\but$,
thus $\bu = \bur + \but$, whereby
\begin{equation}\label{radial and tangential}
 \bur=\stdsum \Ulm \Plm\also \but=\stdsum \Vlm \Blm + \Wlm \Clm.
\end{equation}
We use $\dirrp$ for the vector Dirac delta function on the sphere. Accordingly,
\begin{equation}\label{Dirac delta}
 \intOmega \dirrp\cdot \bu(\rvec)\dOmega=\bu(\rvec').
\end{equation}
The vector spherical-harmonic representation of~$\dirrp$ is the dyad
\begin{equation}\label{Dirac delta expanded}
\dirrp= \stdsum \Plm(\rvec)\Plm(\rvec')+\Blm(\rvec)\Blm(\rvec') +
\Clm(\rvec)\Clm(\rvec').
\end{equation}

\subsection{Bandlimited and Spacelimited Vector Fields}

We consider two subspaces of the space of all square-integrable vector
fields on the unit sphere~$\Omega$. Given
$\sg=(\ldots,\Ulm,\ldots,\Vlm,\ldots,\Wlm,\ldots)\T$, we define the
space of all bandlimited vector fields
$\bandlim=\{\bg:\Ulm=\Vlm=\Wlm=0$ for $L<l\leq \infty$ and $-l\leq
m\leq l$\}, with no power beyond the bandwidth~$L$, whose elements are
the functions
\begin{equation}\label{bandlimited field}
 \bg = \bgr + \bgt= \stdsumL \Ulm\Plm + \Vlm\Blm + \Wlm\Clm
,
\end{equation}
where now
\begin{equation}
 \Ulm=\intOmega \Plm \cdot \bg\dOmega,\label{gUlmVlmWlm}\quad
 \Vlm=\intOmega \Blm \cdot \bg\dOmega,\quad\mbox{and}\,\,
 \Wlm=\intOmega \Clm \cdot \bg\dOmega.
\end{equation}
Similarly, we define $\spacelim=\{\bh: \bh=\ovec$ in~$\Omega\setminus
R \}$ to be the space of all spacelimited vector fields $\bh(\rvec)$
that are equal to zero outside a nonempty region~$R\subset \Omega$. By
definition, the space $\spacelim$ is infinite-dimensional but $\dim
\bandlim = 3(L+1)^2-2$, because the coefficient vector~$\sg$ has
$\sum_{l=0}^L (2l+1)=(L+1)^2$ entries for the $\Ulm$ and $\sum_{l=1}^L
(2l+1)=(L+1)^2-1$ entries for the $\Vlm$ and $\Wlm$, respectively.  We
define the spatial and spectral measures analogously to~(\ref{norms})
\begin{equation}\label{subnorms} \lVert \bg\rVert^2_R = \int_R
  \bg\cdot \bg\dOmega, \qquad \lVert \sg\rVert^2_L= \stdsumL
  \Ulm^2+\Vlm^2 + \Wlm^2.
\end{equation}

\section{Concentration within an Arbitrarily Shaped Region}

No vector field can be strictly bandlimited and strictly spacelimited,
i.e., no $\bu(\rvec)$ can simultaneously be contained in both spaces
$\spacelim$ and $\bandlim$. Our goal is to determine bandlimited
vector fields $\bg(\rvec)\in \bandlim$ with optimal
energy-concentration within a spatial region~$R$, and those
spacelimited vector fields $\bh(\rvec)\in \spacelim$ with an optimally
concentrated spectrum within an interval $0\leq l\leq L$.  Similar to
the scalar time-frequency~\cite{Slepian1960}, multidimensional
Cartesian~\cite{Simons2010,Slepian64} and
spherical~\cite{Simons+2006a} cases, these two spatiospectral
concentration problems are closely related.

\subsection{Spatial Concentration of a Bandlimited Vector  Field}
\label{spatial concentration section} 
 
We maximize the spatial concentration of a bandlimited vector field
$\bg(\rvec)\in \bandlim$ within~$R$ via the ratio
\begin{equation}\label{spatial concentration equation}
 \lambda = \frac{ \lVert \bg\rVert^2_R}{ \lVert
   \bg\rVert^2_\Omega}=\frac{\displaystyle
   \intR\bg\cdot\bg\dOmega}{\displaystyle \intOmega
   \bg\cdot\bg\dOmega }=\text{maximum}. 
\end{equation}
The variational problem~(\ref{spatial concentration equation}) is
analogous to that encountered in one and two scalar dimensions. As
there, the ratio $0<\lambda<1$ is a measure of the spatial
concentration.   

\subsubsection{Purely Radial Vector Fields}
\label{purely radial section}

As a first step, we focus on solving(\ref{spatial concentration
  equation}) for purely radial fields, that is, bandlimited vector
fields in the decomposition~(\ref{radial and
  tangential}),
\begin{equation}\label{bandlg}
\bgr=\stdsumL \Ulm\Plm.
\end{equation}
To simplify the notation we drop the superscript on the
coefficient vector, such that
$\sg=(\ldots,\Ulm,\ldots)\T$ in this Section. 
 Inserting the representation~(\ref{bandlg}) into~(\ref{spatial
  concentration equation}) and switching the order of summation 
and integration, we can express~$\lambda$ as
\begin{equation}
 \lambda = \frac{\displaystyle \stdsumL \Ulm  \stdsumLp \Plmlmp
   \hsom\Ulmp}{\displaystyle \stdsumL  \Ulm^2 }.
\end{equation}
Here we have used orthonormality (\ref{orthonormality})
and the quantities 
\begin{equation}\label{Plmlmp elements}
 \Plmlmp = \intR \Plm\cdot\Plmp\dOmega=\intR \Ylm \Ylpmp\dOmega.
\end{equation}
We can reformulate (\ref{spatial concentration equation}) as a matrix
variational problem~\cite{Horn1990}:
\begin{equation}\label{radial coefficient optimization problem}
 \lambda = \frac{\sg^\sT\sP\hsom\sg}{\sg^\sT\sg} = \text{maximum}
\end{equation}
using the $(L+1)^2\times(L+1)^2$ matrix
\begin{equation}\label{P matrix}
 \sP=\begin{pmatrix}P_{00,00}&\cdots&P_{00,LL}\\
             \vdots&&\vdots\\
	      P_{LL,00}&\cdots&P_{LL,LL}
            \end{pmatrix},
\end{equation}
The stationary solutions of the Rayleigh quotient $\lambda$ in
(\ref{radial coefficient optimization problem}) are solutions of the
$(L+1)^2\times (L+1)^2$ algebraic eigenvalue problem
\begin{equation}\label{radial eigenvalue problem}
 \sP\hsom\sg=\lambda\hsom\sg.
\end{equation}
Therefore the spatial concentration problem of purely radial
bandlimited vector fields is completely equivalent to the scalar
spherical concentration problem~\cite{Simons+2006a}.

\subsubsection{General Vector Fields}

For bandlimited vector fields that are of the kind~(\ref{bandlimited
  field}) described by the complete coefficient  
vector~$\sg = (\ldots,\Ulm,\ldots,\Vlm,\ldots,\Wlm,\ldots)\T$,
operations analogous to those carried out in Section~\ref{purely
  radial section} transform~(\ref{spatial 
  concentration equation}) into a matrix variational problem
in the space of $[3(L+1)^2-2]$-tuples:
\begin{equation}\label{matrix Slepian problem}
 \lambda = \frac{\sg^\sT\sK\hsom\sg}{\sg^\sT\sg} = \text{maximum}.
\end{equation}
Since the inner products of $\Plm$ with $\Blm$
and $\Clm$ are always zero because of~(\ref{strong orthogonality}),
\begin{equation}\label{K matrix}
 \sK=\begin{pmatrix}\sP&\so&\so\\\so&\sB&\sD\\\so&\sD^\sT&\sC\end{pmatrix}=
\begin{pmatrix}\sP&\so\\\so&\sQ\end{pmatrix}, 
\end{equation}
where the $[(L+1)^2-1]\times[(L+1)^2-1]$-dimensional matrices
\begin{equation}\label{B matrix}
 \sB=\begin{pmatrix}B_{10,10}&\cdots&B_{10,LL}\\
             \vdots&&\vdots\\
	      B_{LL,10}&\cdots&B_{LL,LL}
            \end{pmatrix},
\end{equation}
\begin{equation}\label{C matrix}
 \sC=\begin{pmatrix}C_{10,10}&\cdots&C_{10,LL}\\
             \vdots&&\vdots\\
	      C_{LL,10}&\cdots&C_{LL,LL}
            \end{pmatrix},
\end{equation}
\begin{equation}\label{D matrix}
 \sD=\begin{pmatrix}D_{10,10}&\cdots&D_{10,LL}\\
             \vdots&&\vdots\\
	      D_{LL,10}&\cdots&D_{LL,LL}
            \end{pmatrix},
\end{equation}
have matrix entries defined by 
\begin{align}
\Blmlmp=&\intR \Blm\cdot\Blmp\dOmega,\label{Blmlmp elements}\\
\Clmlmp=&\intR \Clm\cdot\Clmp\dOmega,\label{Clmlmp elements}\\
\Dlmlmp=&\intR \Blm\cdot\Clmp\dOmega,\label{Dlmlmp elements}
\end{align}
and
\begin{equation}\label{Q matrix}
\sQ=\begin{pmatrix}\sB&\sD\\\sD^\sT&\sC\end{pmatrix}
=\begin{pmatrix}\hspace{0.65em}\sB&\sD\\-\sD&\sB\end{pmatrix}.
\end{equation}
The last identity follows from (\ref{Blm})--(\ref{Clm}) and
(\ref{Blmlmp elements})--(\ref{Dlmlmp elements}). The solutions to the
concentration problem of general bandlimited vector fields to
arbitrary domains solve the
$[3(L+1)^2-2]\times[3(L+1)^2-2]$-dimensional algebraic eigenvalue
problem
\begin{equation}\label{matrix eigenvalue problem}
 \mathsf{K}\hsom\mathsf{g}=\lambda\hsom\mathsf{g}. 
\end{equation}
The matrix~(\ref{K matrix}) is real, symmetric ($\sK^\sT=\sK$), and it
is positive definite ($\sg^\sT\sK\hsom\sg>0$ for all $\sg\neq \so$),
hence the $3(L+1)^2-2$ eigenvalues~$\lambda$ and associated
eigenvectors~$\sg$ are always real. The eigenvalues
$\lambda_1,\lambda_2,\ldots,\lambda_{3(L+1)^2-2}$ and eigenvectors
$\sg_1,\sg_2,\ldots,\sg_{3(L+1)^2-2}$ can be ordered so that they are
sorted $1>\lambda_1\geq\lambda_2\geq\cdots\geq\lambda_{3(L+1)^2-2}>0$. Every
spectral-domain eigenvector~$\sga$ is associated with a bandlimited
spatial eigenfield~$\bga(\rvec)$ defined by (\ref{bandlimited field}).
If $R$ is a true subset of $\Omega$, then the largest eigenvalue,
$\lambda_1$, will be strictly smaller than one since no bandlimited
function can be non-zero only within a region~$R$ that is smaller than
$\Omega$.  Due to the positive definiteness of the matrix~$\sK$ for a
non-empty region~$R$, the smallest eigenvalue, $\lambda_{3(L+1)^2-2}$,
is larger than zero.

The eigenvectors $\sg_1,\sg_2,\ldots,\sg_{3(L+1)^2-2}$ are
orthogonal. We orthonormalize as
\begin{equation}\label{spectral inner product}
\sga^\sT\hsom\sgb^{}=\deltaab, \qquad 
\sga^\sT\sK\hsom\sgb^{}=\lambda_\alpha\deltaab.
\end{equation}
The associated eigenfields
$\bg_1(\rvec),\bg_2(\rvec),\dots,\bg_{3(L+1)^2-2}(\rvec)$ are a basis
for $\bandlim$ that is orthogonal over the region~$R$ 
and orthonormal over the whole sphere~$\Omega$: 
\begin{equation}\label{spatial inner products}
\intOmega \bga\cdot \bgb \dOmega=\deltaab,\qquad \intR \bga\cdot
\bgb \dOmega = \lambda_\alpha \deltaab. 
\end{equation}
The relations (\ref{spatial inner products}) for the spatial-domain
are equivalent to their matrix counterparts (\ref{spectral inner
  product}).  The eigenfield $\bg_1(\rvec)$ with the largest
eigenvalue~$\lambda_1$ is the element in the space~$\bandlim$ of
bandlimited vector fields with most of its spatial energy within
region~$R$; the eigenfield $\bg_2(\rvec)$ is the next
best-concentrated element in~$\bandlim$ that is orthogonal
to~$\bg_1(\rvec)$ over both~$\Omega$ and~$R$; and so on.

When put into index notation, the eigenvalue
equations~(\ref{matrix eigenvalue problem}) are  
\begin{align}
 \stdsumLp \Plmlmp \Ulmp&=\lambda \Ulm,\label{radial matrix} \\
 \stdsumLp \Blmlmp \Vlmp + \Dlmlmp \Wlmp&=\lambda
 \Vlm,\label{tangential matrix V}\\ 
 \stdsumLp \Dlmlmp^\iT \Vlmp + \Clmlmp \Wlmp&=\lambda
 \Wlm.\label{tangential matrix W} 
\end{align}
By tensor-multiplying the expression~(\ref{radial matrix}) with
$\Plm(\rvec)$, (\ref{tangential matrix V}) with $\Blm(\rvec)$, and (\ref{tangential
  matrix W}) with $\Clm(\rvec)$, and summing in each equation over all
$0\leq l \leq L$ and $-l\leq m \leq l$, we obtain the following system
of spatial-domain equations
\begin{align}
 \intR \left[ \stdsumL \Plm(\rvec) \Plm(\rvec')\right]\cdot
 \bg^r(\rvec')\dOmega' &= \lambda\hsom \bg^r(\rvec),\label{radial
   Fredholm}\\ 
 \intR\left[\stdsumL \Blm(\rvec) \Blm(\rvec')\right]\cdot
 \bg^t(\rvec')\dOmega' &= \lambda \stdsumL
 \Vlm\Blm(\rvec),\label{tangential Fredholm V}\\ 
 \intR\left[ \stdsumL \Clm(\rvec) \Clm(\rvec')\right]\cdot
 \bg^t(\rvec')\dOmega' &= \lambda
 \stdsumL\Wlm\Clm(\rvec).\label{tangential Fredholm W} 
\end{align}
By adding equations~(\ref{radial Fredholm})--(\ref{tangential Fredholm
  W}), we obtain the spatial-domain eigenvalue problem
\begin{equation}\label{Fredholm R}
 \intR \Dirrp\cdot  \bg(\rvec') \dOmega = \lambda\hsom\bg(\rvec), \quad
 \rvec\in \Omega,
\end{equation}
a homogeneous Fredholm integral equation~\cite{Tricomi1970} with a
finite-rank, symmetric, separable, bandlimited vector Dirac delta
function kernel, 
\begin{equation}\label{bandlimited Dirac} 
\Dirrp=\stdsumL\Plm(\rvec)\Plm(\rvec')+\Blm(\rvec)\Blm(\rvec')+
\Clm(\rvec)\Clm(\rvec'), 
\end{equation}
a reproducing kernel~\cite{Simons2010,Simons+2011a} in the space~$\bandlim$. By
inserting the representations (\ref{bandlimited field}) and
(\ref{bandlimited Dirac}) into (\ref{Fredholm R}), we obtain again the
matrix equation (\ref{matrix eigenvalue problem}).
Therefore the
spectral-domain eigenvalue problem for $\sg$ and the spatial-domain
eigenvalue problem for a bandlimited $\bg(\rvec)\in \bandlim$ are
completely equivalent. 

In summary, it is possible to construct an orthogonal family of
bandlimited eigenfields that is optimally concentrated within a
region~$R$ on the unit sphere~$\Omega$ by solving either the Fredholm
integral equation~(\ref{Fredholm R}) for the associated spatial-domain
eigenfields $\bg_1,\bg_2,\ldots,\bg_{3(L+1)^2-2}$ or the matrix
eigenvalue problem~(\ref{matrix eigenvalue problem}) for
spectral-domain eigenvectors $\sg_1,\sg_2,\ldots,\sg_{3(L+1)^2-2}$.
Both methods determine the optimally concentrated eigenfields over the
complete domain $\Omega$, i.e., both in the region~$R$, within which
they are concentrated, and in the complementary
region~$\Omega\setminus R$, into which they show inevitable leakage.

\subsection{Spectral Concentration of a Spacelimited Vector Field}
\label{spectral concentration}

Instead of energy-concentrating a bandlimited vector
field~$\bg(\rvec)\in \bandlim$, within a spatial region~$R$, we may
choose to construct a spacelimited vector field~$\bh(\rvec)\in
\spacelim$ that is concentrated within the spectral interval~$0\leq
l\leq L<\infty$. Such a spacelimited vector function will be defined
by
\begin{equation}
\bh=\stdsum \Ulm'\Plm+\Vlm'\Blm+\Wlm'\Clm,
\end{equation}
with the expansion coefficients given by the spatially limited integrals
\begin{equation}
 \Ulm'=\intR \Plm \cdot \bh\dOmega,\label{hUlmVlmWlm}\quad
 \Vlm'=\intR \Blm \cdot \bh\dOmega,\quad\mbox{and}\,\,
 \Wlm'=\intR \Clm \cdot \bh\dOmega.
\end{equation}

The quadratic concentration measure analogous to~(\ref{spatial
  concentration equation}) is now the ratio 
\begin{equation}\label{spectral concentration equation}
 \lambda=\frac{\lVert \sh\lVert_L^2}{\lVert
   \sh\lVert_\infty^2}=\frac{\displaystyle \stdsumL \Ulm'^2+\Vlm'^2
   +\Wlm'^2}{\displaystyle \stdsum \Ulm'^2+\Vlm'^2
   +\Wlm'^2}=\text{maximum}. 
\end{equation}

The variational problem~(\ref{spectral concentration equation}) can once again be rewritten as a Rayleigh quotient
\begin{equation}
  \lambda = \frac{\displaystyle \intR \intR \bh(\rvec)\cdot \Dirrp\cdot
    \bh(\rvec')\dOmega\dOmega'}{\displaystyle \intR
    \bh(\rvec)\cdot\bh(\rvec)\dOmega}=\text{maximum}
\end{equation}
by inserting the vector spherical-harmonic expansion coefficients (\ref{hUlmVlmWlm}) into (\ref{spectral concentration equation}), switching the order of summation and integration, and by making use of the reproducing property~(\ref{Dirac delta})
of the delta function (\ref{Dirac delta expanded}) and the definition
(\ref{bandlimited Dirac}) of the kernel $\Dirrp$.
Stationary solutions of~(\ref{spectral concentration
  equation}) solve the Fredholm integral
equation
\begin{equation}\label{Fredholm limited space}
 \intR \Dirrp\cdot \bh(\rvec')\dOmega' = \lambda\hsom\bh(\rvec), \quad \rvec\in R.
\end{equation}
This equation for $\bh(\rvec)\in \spacelim$ is identical to
(\ref{Fredholm R}) for $\bg(\rvec)\in \bandlim$, the difference being
that~(\ref{Fredholm R}) is applicable on the entire sphere~$\Omega$,
while the domain of~(\ref{Fredholm limited space}) is limited to the
region~$R$, within which $\bh(\rvec)\neq\ovec$.  We constructed the
spectral norm ratio maximizing eigenfields $\bh(\rvec)$ for
(\ref{spectral concentration equation}) such that they are identical
to the eigenfields $\bg(\rvec)$ that maximize the spatial norm ratio
(\ref{spatial concentration equation}) within the region~$R$.  We
normalize such that
\begin{equation}\label{h normalization}
 \bh(\rvec)=\begin{cases} \bg(\rvec) &\text{if } \rvec\in R,\\
				    \ovec              &\text{otherwise}.
		      \end{cases}
\end{equation}
Every bandlimited eigenfield~$\bga\in \bandlim$ leads to a
spacelimited~$\bha\in \spacelim$ by the restriction (\ref{h
  normalization}).  The eigenvalues $\lambda_\alpha$ associated with
the corresponding~$\bga$ measure the fractional spatial energy
$1-\lambda_\alpha$ that leaked to the region~$\Omega\setminus R$.
These eigenvalues are identical to the fractional spectral energy that
leaked into the degrees $L<l\leq \infty$ by truncating $\bga$ in the
construction of $\bha$ (\ref{h normalization}).  Equivalently, we
could have started with the variational problem~(\ref{spectral
  concentration equation}) instead of~(\ref{spatial concentration
  equation}) to obtain the integral equation~(\ref{Fredholm R}) and
then extend the domain (\ref{Fredholm limited space}) to the whole
sphere~$\Omega$.

The constructed spacelimited eigenfields
$\bh_1(\rvec),\bh_2(\rvec),\ldots,\bh_{3(L+1)^2-2}$ from (\ref{h
  normalization}) are orthogonal over both the
whole sphere~$\Omega$ and the region~$R$: 
\begin{equation}
 \intOmega \bha\cdot \bhb \dOmega = \intR \bha \cdot\bhb \dOmega =
 \lambda_\alpha\deltaab. 
\end{equation}
We can express the coefficients of
$\sh=(\ldots,\Ulm',\ldots,\Vlm',\ldots,\Wlm',\ldots)\T$, where $0\leq
l\leq \infty$ by the coefficients
$\sg=(\ldots,\Ulm,\ldots,\Vlm,\ldots,\Wlm,\ldots)\T$, with $0\leq
l\leq L$ using the relation $\sh=\sK\hsom\sg$, which leads to
$\Ulm'=\lambda\Ulm$, $\Vlm'=\lambda\Vlm$ and $\Wlm'=\lambda\Wlm$, when
$0\leq l\leq L$ due to (\ref{radial matrix})--(\ref{tangential matrix
  W}).  The solutions to equation~(\ref{Fredholm limited space}) form
an infinite-dimensional space. The complement to the $3(L+1)^2-2$
eigenfields with nonzero eigenvalues
$\lambda_1,\lambda_2,\ldots,\lambda_{3(L+1)^2-2}$ is the space spanned
by all eigenfields of~(\ref{Fredholm limited space}) with associated
eigenvalue $\lambda=0$.  Fields $\bh(\rvec)$ vanishing
in~$\Omega\setminus R$ without power in the spectral interval $0\leq
l\leq L$ are members of this null space.

\subsection{Significant and Insignificant Eigenvalues}

The eigenvalues of the matrix $\sK$ defined in (\ref{K matrix}) can be
summed up as follows
\begin{align}\label{Shannon number}
 N&=\sumalpha \lambda_\alpha = \tr \sK= \stdsumL
  (\Plmlm+\Blmlm+\Clmlm) \\ &= \intR \left[\stdsumL
  \Plm(\rvec)\cdot \Plm(\rvec)
  +\Blm(\rvec)\cdot \Blm(\rvec) +\Clm(\rvec)\cdot \Clm(\rvec)\right] \dOmega\\
  &=\left[3(L+1)^2-2\right]\frac{A}{4\pi}.
\end{align}
In the fourth equality we substituted the diagonal matrix elements $\Plmlm,
\Blmlm$ and $\Clmlm$ from (\ref{Plmlmp elements}), (\ref{Blmlmp
  elements})--(\ref{Clmlmp elements}), and in the last equality we
used the addition theorems~(\ref{vecadd1})--(\ref{vecadd3}).

The value~$N$ in~(\ref{Shannon number}) is the vector spherical
analogue of the Shannon number in the scalar Slepian concentration
problems~\cite{Simons2010}.  Well conctentrated eigenfields
$\bga(\rvec)$ for the region~$R$ will have
eigenvalues~$\lambda_\alpha$ near unity, whereas poorly concentrated
eigenfields will have eigenvalues $\lambda_\alpha$ close to zero. Due
to the characteristic step-shaped spectrum of eigenvalues\\
$\lambda_1,\lambda_2,\ldots,\lambda_{3(L+1)^2-2}$, the total number of
significant ($\lambda_\alpha \approx 1$) eigenvalues can be well
approximated by the rounded sum~(\ref{Shannon number}), as in the
one-dimensional and two-dimensional scalar spherical problems.
Since~$N$ is a good estimate for the number of significant
eigenvalues, then, roughly speaking, the vector spherical Shannon
number~(\ref{Shannon number}) describes the dimension of the space of
vector fields~$\bu(\rvec)$ that are approximately limited in both the
spectral domain to vector spherical-harmonic degrees $0\leq l\leq L$,
and in the spatial domain to an arbitrarily shaped region~$R$ of
area~$A$~\cite{Landau1965,Landau1967}.

Instead of constructing a bandlimited field $\bg(\rvec)\in \bandlim$
that is optimally energy-concentrated within a spatial region~$R$, we
could have sought to construct one that is optimally excluded
from~$R$, i.e., one that is optimally concentrated
within~$\Omega\setminus R$ and therefore sought to minimize rather
than maximize the Rayleigh quotient~(\ref{spatial concentration
  equation}).  What we have constructed are the stationary solutions
$\bg(\rvec)\in \bandlim$ of (\ref{spatial concentration equation}).
Therefore we have actually solved the concentration and exclusion
problems simultaneously. The optimally excluded eigenfields are
identical to the optimally concentrated eigenfields but with reversed
ordering.  Because $\lambda_\alpha$ is the fractional power of $\bga$
within~$R$, its fractional power within $\Omega \setminus R$ is
$1-\lambda_\alpha$ . If the region~$R$ of area~$A$ covers only a small
fraction of the sphere $A\ll4\pi$, the number of well-excluded
eigenfields will be much larger than the number of well concentrated
eigenfields.

We can express 
the kernel~$\Dirrp$ in the integral eigenvalue equation~(\ref{Fredholm
  R}) in terms of the spatial-domain eigenfields
$\bg_1,\bg_2,\ldots,\bg_{3(L+1)^2-2}$ in the form 
\begin{equation}\label{Mercer}
 \Dirrp = \sumalpha \bga(\rvec)\bga(\rvec').
\end{equation}
Equation~(\ref{Mercer}) is equivalent to the original
representation~(\ref{bandlimited Dirac}), because both the $\Plm,
\Blm, \Clm$, $0\leq l\leq L$, $-l\leq m \leq l$, and the
$\bga,\alpha=1,2,\ldots,3(L+1)^2-2$, are $[3(L+1)^2-2]$-dimensional
orthonormal bases for $\bandlim$, and the transformation matrix that
consists of the eigenvectors is orthogonal. The transformed
representation~(\ref{Mercer}) is a vector spherical version of
Mercer's theorem~\cite{Tricomi1970,Flandrin1999,Kanwal1971}.  Upon
setting $\rvec' = \rvec$ in~(\ref{Mercer}) and applying the trace
  \cite{Dahlen1998}, we deduce that the sum of the squares of the
  $3(L+1)^2-2$ bandlimited eigenfields $\bga(\rvec)$ is a constant that is
  independent of position $\rvec$ on the sphere~$\Omega$,
\begin{equation}
  \sumalpha
  \bga(\rvec)\cdot\bga(\rvec)=\frac{3(L+1)^2-2}{4\pi}=\frac{N}{A}. 
\end{equation}
If the eigenvalues of the first $N$ eigenfields
$\bg_1(\rvec),\bg_2(\rvec),\ldots,\bg_N(\rvec)$ are near
unity, and the remaining eigenvalues
$\bg_{N+1}(\rvec),\bg_{N+2}(\rvec),\ldots,\bg_{3(L+1)^2-2}(\rvec)$
are near zero,
then we expect the eigenvalue-weighted sum of squares to be 
\begin{equation}\label{cumulative energy}
 \sumalpha \lambda_\alpha\, \bga(\rvec)\cdot\bga(\rvec)\approx
 \sum_{\alpha=1}^N \lambda_\alpha\,\bga(\rvec)\cdot\bga(\rvec) 
\approx\begin{cases}
    N/A&\text{if }\rvec\in R,\\
    0  &\text{otherwise}.
\end{cases}
\end{equation}
The terms with $N+1\leq \alpha\leq 3(L+1)^2-2$ should be comparatively
small.  It is hence immaterial whether we include them in the sum
(\ref{cumulative energy}) or not. The combination of the first $N$
orthogonal eigenfields $\bga,\alpha=1,2,\ldots,N$, with eigenvalues
$\lambda_\alpha\approx 1$, provide an essentially uniform coverage of
the region~$R$.  This characterizes the spatiospectral concentration
problem. The spatiospectrally concentrated basis reduces the number of
degrees of freedom from $\dim \bandlim = 3(L+1)^2-2$ to
$N=[3(L+1)^2-2]A/(4\pi)$.

\subsection{Pairs of Spatially Concentrated  Tangential Vector Fields}
\label{eigenvalue pairs}

It is possible to construct, from one spatially concentrated,
bandlimited tangential vector field another orthogonal, equally
concentrated and equally bandlimited vector field, by simply rotating
its vectorial directions at each point on the sphere by $90^\circ$
while retaining the absolute values. Such pairs of tangential Slepian
fields already appear in the purely tangential eigenvalue problem,
which, due to the block-diagonal shape of~$\sK$ in~(\ref{K matrix}),
can be solved independently from the radial problem~(\ref{radial
  eigenvalue problem}). As can be seen from (\ref{Plm})--(\ref{Clm})
and (\ref{B matrix})--(\ref{Dlmlmp elements}), the following holds for
any region~$R$ and bandlimit~$L$: 
\begin{equation}\label{knowalso}
\sC = \sB \text{ and } \sD^\sT = -\sD,
\end{equation}
and thus the purely tangential concentration problem is
\begin{equation}\label{tangential matrix}
\sQ\hsom\sg=\begin{pmatrix}{\hspace{0.75em}}\sB&\sD\\-\sD&\sB\end{pmatrix}\sg
= \lambda\hsom\sg. 
\end{equation}
If $\sg = (\sg_1, \sg_2)^\sT$ is an eigenvector of~(\ref{tangential
  matrix}) with eigenvalue $\lambda$, then $\sg = (-\sg_2, \sg_1)^\sT$
is also an eigenvector with the same associated eigenvalue $\lambda$.
The Slepian field constructed from $(-\sg_2, \sg_1)^\sT$ has the same
pointwise absolute value as the Slepian field constructed from
$(\sg_1, \sg_2)^\sT$, and they are pointwise orthogonal.

 
\section{Concentration within an axisymmetric polar cap}
\label{polar cap section}

In this Section we concentrate on the special but important case where
$R$ is a symmetric polar cap with colatitudinal radius~$\Theta$ that
is centered on the north pole, as is shown in Figure~\ref{diagram}.
Because rotations on the sphere commute with the
operators~(\ref{Plm})--(\ref{Clm}) defining the vector spherical
harmonics~\cite{Freeden2009}, the optimally concentrated eigenfields
of the polar cap $R=\{\theta: 0 \leq \theta \leq \Theta\}$ can be
rotated to anywhere on the unit sphere using the same transformations
that apply in the rotation of scalar functions
\cite{Blanco+97,Dahlen+98,Edmonds96}.

\subsection{Decomposition of the Spectral-Domain Eigenvalue Problem}

In the axisymmetric case the matrix elements~(\ref{Plmlmp elements})
and (\ref{Blmlmp  elements})--(\ref{Dlmlmp elements}) reduce to  
\begin{eqnarray}
 \Plmlmp&=&2\pi\hsom\deltamm\intoth \Xlm\Xlpm \sindth,\\
 \Blmlmp&=&\frac{\displaystyle 2\pi\hsom\deltamm\intoth 
\left[\dXlm\dXlpm+m^2\divsinsq \Xlm \Xlpm \right] \sindth}{\displaystyle
   \renormprod},\\ 
\Dlmlmp&=&-\frac{\displaystyle 2\pi\hsom\delta_{-mm'}\hsom m
  \Xlm(\Theta)\Xlpm(\Theta)}{\displaystyle \renormprod} 
,
\end{eqnarray}
while, as we know from~(\ref{knowalso}) also, $\Clmlmp=\Blmlmp$ and we
remember~(\ref{remind}).

The Kronecker deltas $\deltamm$ and $\delta_{-mm'}$ admit rearranging
the ${(L+1)^2\times(L+1)^2}$ radial-compo\-nent  
matrix $\sP$ and the ${[2(L+1)^2-2]} \times[2(L+1)^2-2]$
tangential-compo\-nent matrix $\sQ$ such that they are block-diagonal:
$\sP=\diag(\sP_0,\sP_1, \sP_{-1},\ldots,\sP_L,\sP_{-L})$ and  
$\sQ=\diag(\sQ_0,\sQ_1, \sQ_{-1},\ldots,\sQ_L,\sQ_{-L})$. 

Instead of solving the full eigenvalue equation~(\ref{matrix
eigenvalue problem}), we can thus elect to solve a series of smaller
spectral-domain algebraic eigenvalue problems, one for each order, 
\begin{equation}\label{radial and tangential polar cap eigenvalue
problems} 
\sPm\hsom\sg=\lambda\hsom\sg_m  \also
\sQm\hsom\sg=\lambda\hsom\sg_m.\end{equation}

The matrices $\sPm$ and $\sQm$ 
are of the form 
\begin{equation}
 \sPm=\begin{pmatrix}
   P^m_{mm}&\cdots&P^m_{mL}\\\vdots&&\vdots\\P^m_{Lm}&\cdots&P^m_{LL}
 \end{pmatrix}, \quad  
\sQm= \begin{pmatrix}\sBm&\sDm\\\sDm\T&\sBm\end{pmatrix}, 
\end{equation} 
with
\begin{equation}
\sBm=\begin{pmatrix}
  B^m_{mm}&\cdots&B^m_{mL}\\\vdots&&\vdots\\B^m_{Lm}&\cdots&B^m_{LL}
\end{pmatrix}, \quad 
\sDm=\begin{pmatrix}
  D^m_{mm}&\cdots&D^m_{mL}\\\vdots&&\vdots\\D^m_{Lm}&\cdots&D^m_{LL}
\end{pmatrix},  
\end{equation}
where, for any particular harmonic order $0\leq m \leq L$ and degree
$m\leq l\leq L$, we denote $\Pllp^m=\Plmlpm$, 
and, likewise, for ${\maxm\leq l\leq L}$, we denote
$\Bllp^m=\Blmlpm,$ and $\Dllp^m=\Dlmlpm$. 
We then also have  
\begin{equation}
\sPmm=\sPm, \quad \sBmm=\sBm, \quad \sDmm=-\sDm, \quad \sD_0=\so,
\end{equation}
and that $\sPm,\sBm,\sDm,$ and consequently, $\sQm$, are symmetric.

The calculations of the matrix elements $\Dllp^m$ are straightforward
since they merely consist in evaluating the~$\Xlm$ at~$\Theta$. The
calculations of the product integrals $\Pllp^m$ and $\Bllp^m$ can be
simplified to integrations over individual terms~$\Xlm$. For example,
for the elements $\Pllp^m$ we can directly apply~(\ref{Gaunt}) to
reduce them to 
\begin{align}
\nnr\Pllp^m=& \sqrt{\pi(2l+1)(2l'+1)} \gauntsum \sqrt{2n+1}
\wigner{l}{n}{l'}{0}{0}{0}
\wigner{l}{n}{l'}{m}{-2m}{m}\\ 
&\hspace{12em}\times
 \intoth X_{n\,2m}(\theta)\sin \theta \,d\theta,\label{Gaunt sum
   integration} 
\end{align}
and those can be handled recursively via equation~(\ref{bladipaul}).
Here, as before, we set $\Xlm=0$ for $m>l$.
Since this solution to~(\ref{radial eigenvalue problem}) is identical
to that of the scalar spherical concentration problem for the polar
cap, alternate expressions and special cases can be found
elsewhere~\cite{Simons+2006a,Wieczorek+2005}. For the $\Bllp^m$ at
positive orders $m\geq 0$, as first noted by
Eshagh~\cite{Eshagh2009a}, we first need to transform the derivative
products~$\dXlm \dXlpm$, and $m^2 \divsinsq \Xlm\Xlpm$, into products
of $\Xlm$ using the lemmas~(\ref{Ilk 1})--(\ref{ilk coefficients b}),
to
\begin{eqnarray}\label{XpXp}
\lefteqn{\intoth \dXlm \dXlpm \sindth=}\\
  \nnr&&\hspace{4.5em}
 \amlm \amlpm \intoth X_{l\, m- 1}X_{l'\, m- 1}\sindth
+\aplm \amlpm \intoth X_{l\, m+ 1}X_{l'\, m- 1}\sindth\\
\nnr &&\hspace{4.5em}
{}+\amlm \aplpm \intoth X_{l\, m- 1}X_{l'\, m+ 1}\sindth
+\aplm \aplpm \intoth X_{l\, m+ 1}X_{l'\, m+ 1}\sindth,
\end{eqnarray}
\begin{eqnarray}\label{msinXX}
\lefteqn{\intoth m^2 \divsinsq\, \Xlm \Xlpm \sindth=}\\ 
  \nnr&&\hspace{1em}
\bmlm \bmlpm \intoth X_{l-1\, m- 1}X_{l'-1\, m- 1}\sindth
+\bplm \bmlpm \intoth X_{l-1\, m+ 1}X_{l'-1\, m- 1}\sindth\\
\nnr&&\hspace{1em}
{}+\bmlm \bplpm \intoth X_{l-1\, m- 1}X_{l'-1\, m+ 1}\sindth
+\bplm \bplpm \intoth X_{l-1\, m+ 1}X_{l'-1\, m+ 1}\sindth,
\end{eqnarray}
where $\apmlm$ and $\bpmlm$ are defined in~(\ref{ilk coefficients a})
and~(\ref{ilk coefficients b}). The right hand sides of~(\ref{XpXp})
and~(\ref{msinXX}) can be expanded using~(\ref{Gaunt}) and then the
recursion~(\ref{bladipaul}) can be applied.

We order the $L-m+1$ distinct eigenvalues of $\sPm$ and the
${2(L-\maxm+1)}$ distinct eigenvalues of $\sQm$ obtained by solving
each of the eigenvalue problems~(\ref{radial and
  tangential polar cap eigenvalue problems}) so that
$1\geq\lambda_1\geq\lambda_2\geq \cdots >0$. Additionally we orthonormalize the
associated eigenvectors $\sg_1,\sg_2,\ldots$ as in~(\ref{spectral
  inner product}) so that
\begin{equation}
 \sga^\sT\hsom\sgb = \deltaab, \qquad \sga^\sT\sP\hsom\sgb =
 \lambda_\alpha\deltaab \quad\text{or}\quad   
 \sga^\sT\sQ\hsom\sgb 
 = \lambda_\alpha\deltaab,
\end{equation}
depending on whether $\sga$ and $\sgb$ are the eigenvectors of
$\sPm$ or of $\sQm$.  

\begin{figure}
\begin{centering}
  \includegraphics[width=0.85\textwidth]{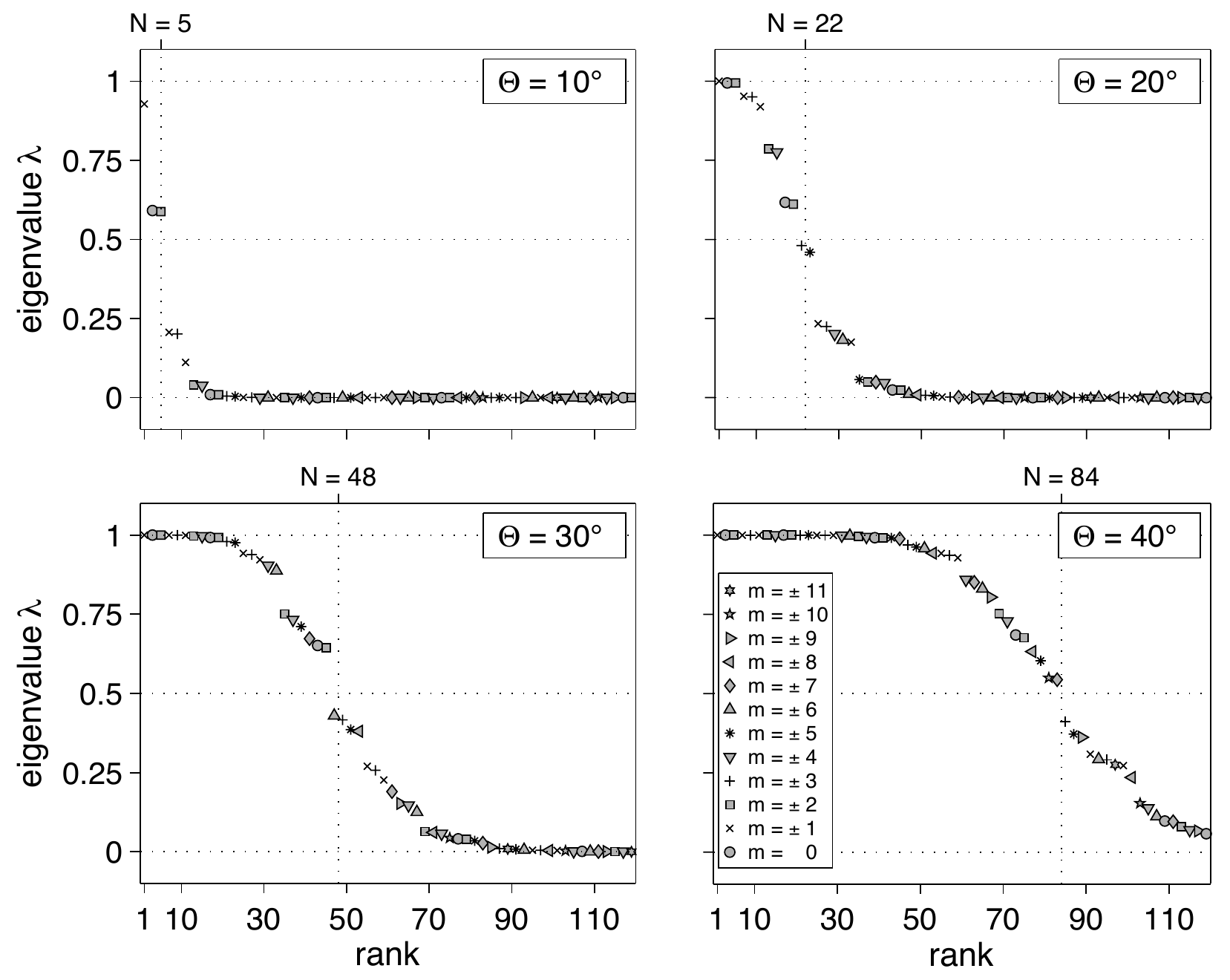}\figend
  \caption{\label{psvals} Reordered eigenvalue spectra
    ($\lambda_\alpha$ versus rank $\alpha$) for the tangential vector
    Slepian functions for axisymmetric polar caps of colatitudinal
    radii $\Theta=10^\circ, 20^\circ, 30^\circ, 40^\circ$ and a common
    bandwidth $L = 18$. The total number of eigenvalues is $2(L+1)^2-2
    = 720$; only $\lambda_1$ through $\lambda_{120}$ are shown.
    Different symbols are used to plot the orders $-11 \leq m \leq
    11$. Each symbol stands for two eigenvalues, that is, the $\pm m$
    doublets for $m>0$ and the doublets stemming from the
    block-diagonal shape of $\sQ_0$ for $m=0$. Vertical gridlines and
    top labels specify the rounded Shannon numbers $N^t = 5, 22, 48$,
    and~$84$.  }
\end{centering}
\end{figure}

\begin{figure}
\begin{centering}
 \includegraphics[height=0.80
 \textheight]{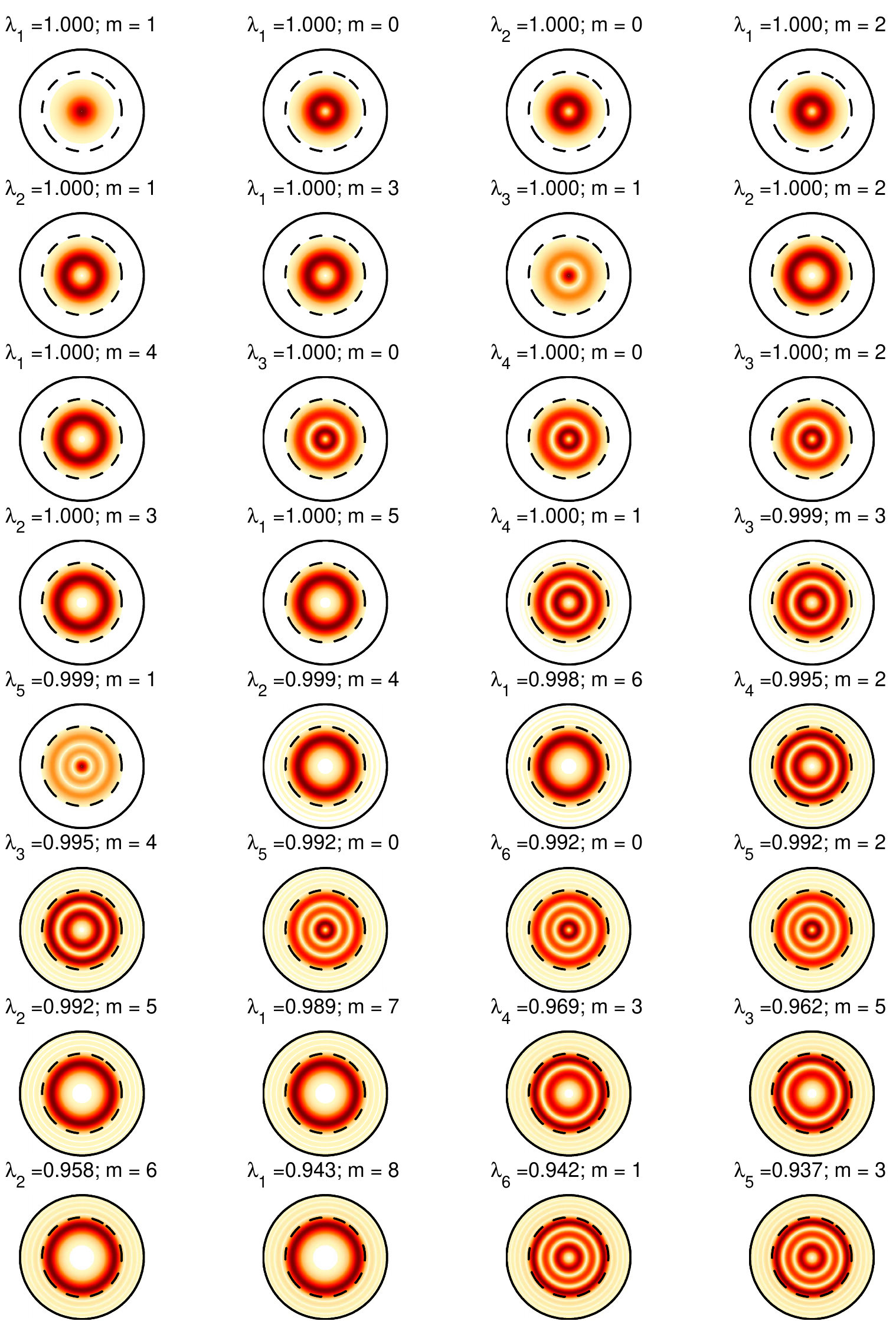}\figend
 \caption{\label{absolute polar cap} Absolute values of the tangential bandlimited 
 eigenfields $\lvert \bg(\theta,\phi)\rvert$ that are optimally concentrated within a circular
 cap of colatitudinal radius $\Theta = 40^\circ$. Dashed circles denote the cap boundary.
 The bandwidth is $L=18$, and the rounded Shannon number
   for the tangential space $N^t=84$. Subscripts on the eigenvalues $\lambda_\alpha$ 
   specify the fixed-order rank.  Only absolute values
   for $m\geq 0$ are shown because the absolute values for $\pm m$ are
   identical. The eigenvalues have been resorted into a mixed-order ranking, with the 
   best-concentrated eigenfields plotted on the top left and a decreasing concentration ratio
   to the right and downwards. Regions in which the absolute value is less than one
   hundredth of the maximum value on the sphere are left white.}
\end{centering}
\end{figure}

\subsection{Eigenvalue Spectrum and Eigenfields}
\label{eigspec}

For the  fixed-order radial eigenvalue problem~(\ref{radial and
  tangential polar cap eigenvalue problems}) we can calculate the number of
significant eigenvalues, or partial Shannon number,
using any of the two formulas 
\begin{equation}
N^r_m=\sum_{\alpha=1}^{L-m+1}\lambda_\alpha=\sum_{l=m}^L \Pll^m. 
\end{equation}
For the fixed-order tangential eigenvalue problem~(\ref{radial and
  tangential polar cap eigenvalue problems}) we obtain the number
of significant eigenvalues from
\begin{equation}
N^t_m=\sum_{\alpha=1}^{2(L-\maxm+1)}\lambda_\alpha=2\sum_{l=\maxm}^L \Bll^m.
\end{equation}
Once we have found the $L+1$ sequences of fixed-order radial and
tangential eigenvalues, we can resort them into an overall mixed-order
ranking. The Shannon number of radial eigenvalues is then given by
$N^r=N^r_0+2\sum_{m=1}^L N^r_m$, while that of the tangential
eigenvalues is $N^t=N^t_0+2\sum_{m=1}^L N_m$. The factor of two
accounts for the $\pm m$ degeneracy; the total number of significantly
concentrated vector fields is $N=N^r+N^t$.

Figure~\ref{psvals}, shows the reordered, mixed-$m$ eigenvalue spectra
for four different polar caps, with colatitudinal radii $\Theta =
10^\circ, 20^\circ, 30^\circ, 40^\circ$ in the case of the tangential
concentration problem. The spherical-harmonic bandlimit is $L=18$.
Each symbol stands for two eigenvalues, arising from the
plus-minus-degeneracy for $m>0$ and from the block-diagonal shape for
$m=0$.  The spectra have a characteristic step shape
\cite{Landau1965,Percival1993,Slepian1965}, showing significant
($\lambda\approx 1$) and insignificant ($\lambda\approx 0$)
eigenvalues separated by a transition band. In all four cases the
reasonably well concentrated eigensolutions $(\lambda \geq 0.5)$ and
the more poorly concentrated ones $(\lambda < 0.5)$ are separated by the 
rounded Shannon number.

Figure~\ref{absolute polar cap} shows a polar plot of the absolute
values of the first 32 tangential eigenfields~$\lvert
\bg(\rvec)\rvert$ constructed from the eigenvectors of $\sQ$.  The
reconstructed fields for positive and negative orders $\pm m$ have the
same absolute values. We therefore only plot the absolute values for
$m\geq 0$. The eigenfields are concentrated within a cap of radius
$\Theta = 40^\circ$. The maximal vector spherical-harmonic degree is
$L=18$ and the Shannon number of the tangential problem $N^t=84$. The
eigenvalue ranking is mixed-order and the concentration factors
$1<\lambda\leq 0.937$ and orders $m$ of each absolute field are
indicated. Red denotes the maximum value while all absolute values
smaller than~$1\%$ of the maximum value are white.

\begin{figure}
\begin{centering}
 \includegraphics[width=0.75\textwidth]{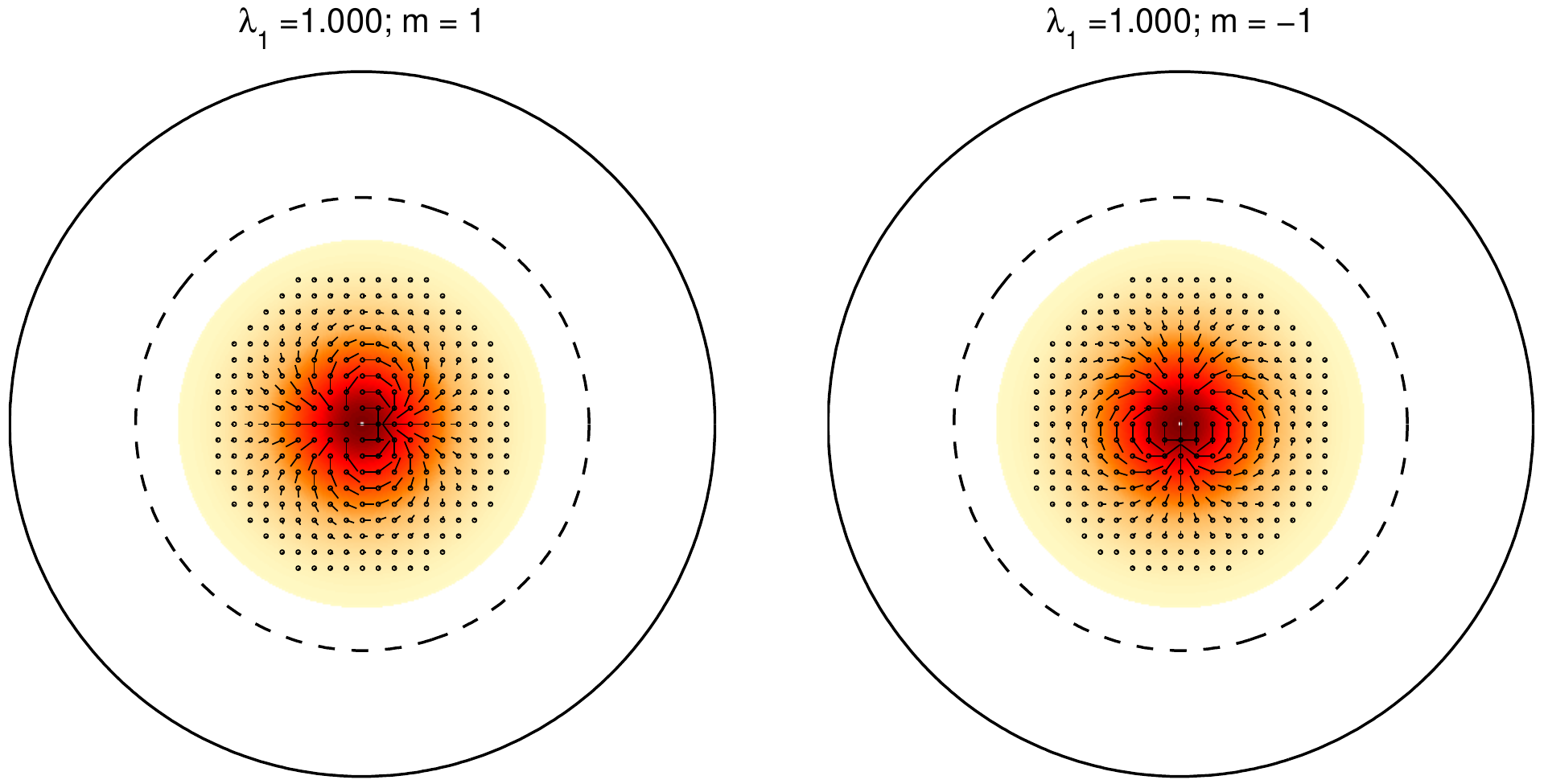}\figend
 \caption{\label{tangential polar cap}Bandlimited tangential Slepian
   functions~$\bg(\theta,\phi)$, of spherical-harmonic orders ${m=\pm
     1}$, optimally concentrated within a polar cap of radius $\Theta
   = 40^\circ$. The bandwidth is $L=18$. Color is absolute value (red
   the maximum) and circles with strokes indicate the direction of the
   eigenfield on the tangential plane. Regions in which the absolute
   value is less than one hundredth of the maximum absolute value on
   the sphere are left white.}
\end{centering}
\end{figure}

In Figure~\ref{tangential polar cap} we illustrate the
best-concentrated tangential Slepian fields for order~$|m|=1$,
corresponding to the overall best-concentrated absolute field of
Figure~\ref{absolute polar cap}.  The left panel shows the vector
field for $m=1$, thus the reconstruction using the best-concentrated
eigenvector of~$\sQ_1$. The right panel shows the vector field for
$m=-1$, thus the reconstruction using the best-concentrated
eigenvector of~$\sQ_{-1}$.  As in Figure~\ref{absolute polar cap}, the
radius of the polar cap is $\Theta=40^\circ$ and the bandwidth $L=18$.
Both vector fields have a singularity at the north pole $\theta=0$.
This is due to the fact that all the $\Blm$ and $\Clm$ for $m=1$ have
a singularity at the north pole stemming from the derivatives of the
$\Xlm$ (\ref{Ilk 1}) which are not equal to zero at $\theta=0$. The
dashed circles denote the cap boundary. The color scales with the
absolute value of the vector field, ranging from white for values
below $1\%$ of the maximum to red for the maximum value. The
directions of the field are indicated by accordingly oriented strokes
at the positions marked by the open circles.

\begin{figure}
\begin{centering}
 \includegraphics[angle=0, width=0.85\textwidth]{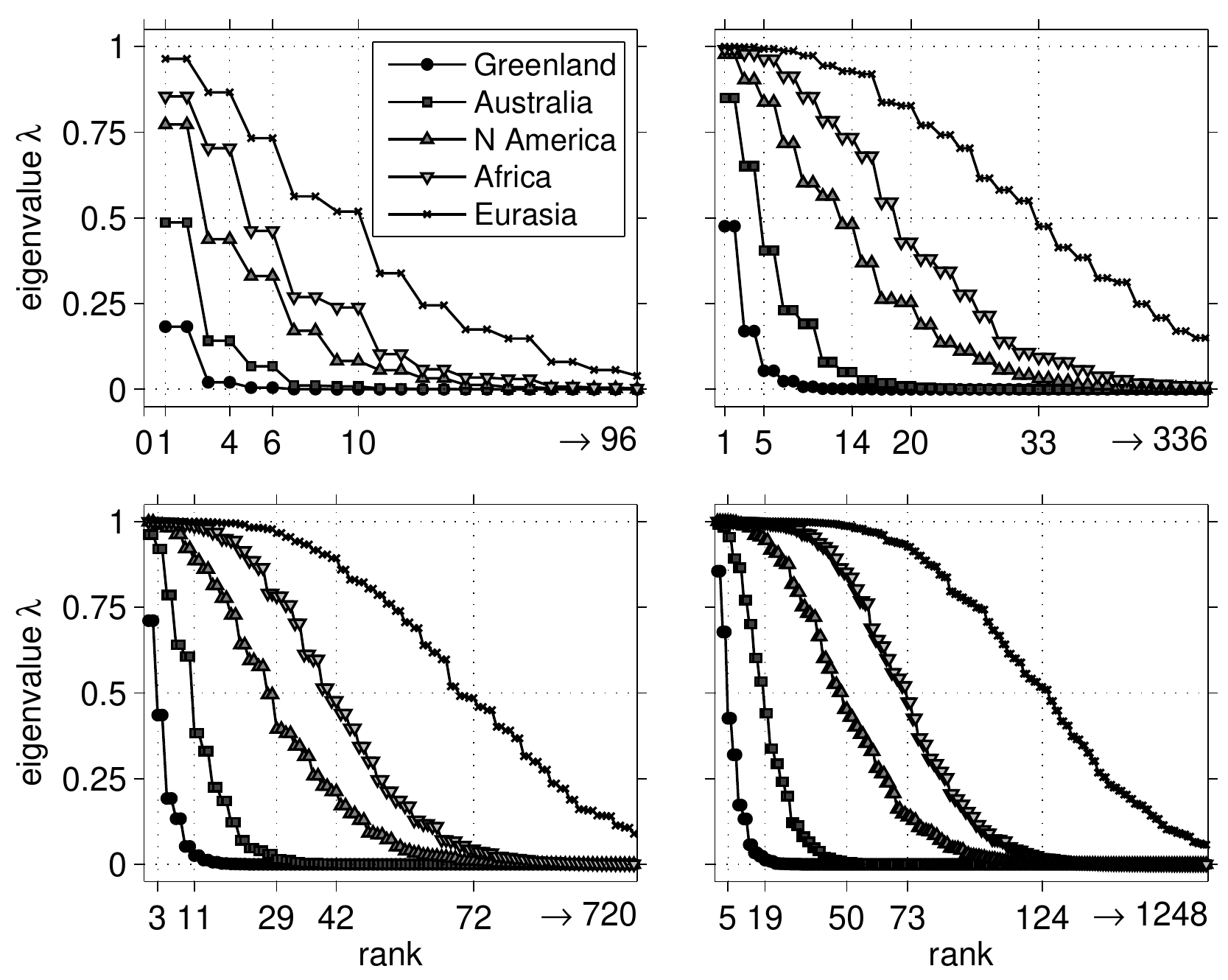}\figend
 \caption{\label{regions eigenvalues}Eigenvalue spectra for the
   tangential fields for Greenland, Australia, North America, Africa,
   and Eurasia. From upper left to lower right, four bandwidths,
   $L=6,12,18,24$, are considered. The horizontal axis in each panel
   is truncated; the total number of eigenvalues
   $2(L+1)^2-2=96,336,720,1248,$ appears to the right of the arrow.
   Vertical gridlines and the five leftmost ordinate labels specify
   the rounded Shannon numbers~$N^t$.}
\end{centering}
\end{figure}

\section{Continental Concentration}

In the following example we consider the spatiospectral concentration
in seven of Earth's continental regions. Together with their rounded
tangential-component Shannon numbers, $N^t=[2(L+1)^2-2]A/(4\pi)$, the
regions are listed in Table~\ref{cont table} for different bandwidths.
The (tangential) spherical Slepian fields that we will be showing
should be well suited to the localized analysis of global
vector-valued satellite-magnetic data such as measured to study the
magnetization of the terrestrial
lithosphere~(e.g.~\cite{Friis2006,Maus2006,Beggan+2013,Sabaka2006}),
or more generally, planetary magnetic
fields~\cite{Schott+2011,Langlais+2010,Lewis+2012,Sterenborg+2007}.

\begin{table}[h]
  \caption{\label{cont table}Fractional areas, tangential Shannon
    numbers, and bandwidths for the vectorial concentration problem to
    continental areas.}
\hrule\vspace{1em}
\begin{center}
\begin{tabular}{l|c|cccc}
\hline
&Fractional area &\multicolumn{4}{c}{ Tangential Shannon number $N^t$ }\\
Continental region&$A/(4\pi)$ in \%&$L=6$&$L=12$&$L=18$&$L=24$\\	
\hline
Greenland&0.43&0&1&3&5\\
Australia&1.51&1&5&11&19\\
Antarctica&2.72&3&9&20&34\\
South America&3.45&3&12&25&43\\
North America&4.03&4&14&29&50\\
Africa&5.81&6&20&42&73\\
Eurasia&9.97&10&33&72&124\\
\hline
\end{tabular}
\end{center}
\end{table}


\subsection{Bandlimited Fields}

Figure~\ref{regions eigenvalues} shows the eigenvalue spectra of the
tangential Slepian fields for the five regions Greenland, Australia,
North America, Africa, and Eurasia, and four spherical-harmonic
bandlimits, $L=6,12,18,24$, which correspond to $2(L+1)^2-2 = 96$,
336, 720, 1248 eigenfields each. The smallest wavelength for a
bandwidth limit~$L$ is $2\pi/\sqrt{L(L+1)}\approx 2\pi/(L+1/2)$
multiplied by Earth's radius~\cite{Jeans1923}. The cutoff wavelengths
for $L=6,12,18$, and $24$ are 6200, 3200, 2200, and 1600~km,
respectively.  Only Eurasia, the largest region, has enough area to
contain at least one nearly perfectly concentrated eigenfield for the
smallest bandwidth, $L=6$, and Greenland, the smallest of the
considered regions, is too small to contain even for the largest
bandwidth, $L=24$, a single eigenfield with a concentration factor
$\lambda$ near unity. Again, as was the case for a polar cap
(Figure~\ref{psvals}), the well-concentrated eigenfunctions with
eigenvalues $\lambda\geq 0.5$ are separated from the poorly
concentrated ones with eigenvalues $\lambda < 0.5$ by the rounded
Shannon numbers~$N^t$.  The eigenvalues occur in pairs as described in
Sections~\ref{eigenvalue pairs} and~\ref{eigspec} 

\begin{figure}[h]
\begin{centering}
 \includegraphics[angle=0, width=0.665\textwidth]{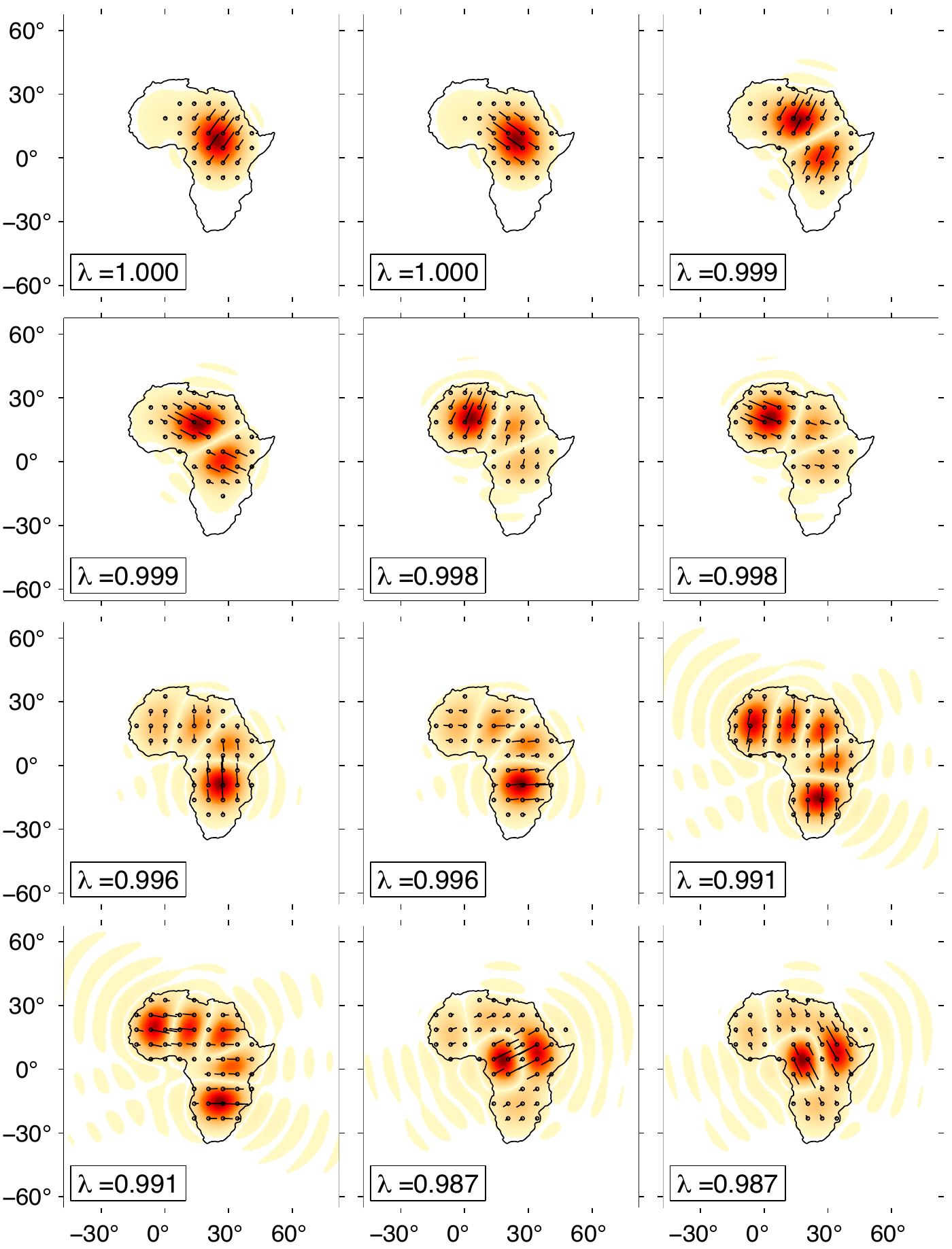}\figend 
 \caption{\label{tangential africa} Twelve tangential Slepian
   functions $\bg_1,\bg_2,\ldots,\bg_{12}$, bandlimited to $L=18$,
   optimally concentrated within Africa. The concentration factors
   $\lambda_1, \lambda_2,\ldots,\lambda_{12}$ are indicated. The
   rounded tangential Shannon number $N^t=42$. Order of concentration is
   left to right, top to bottom. Color scheme and symbols are as in
   Figure~\ref{tangential polar cap}. 
 }
\end{centering}
\end{figure}

\begin{figure}[h]
\begin{centering}
  \includegraphics[angle=0,
  width=0.675\textwidth]{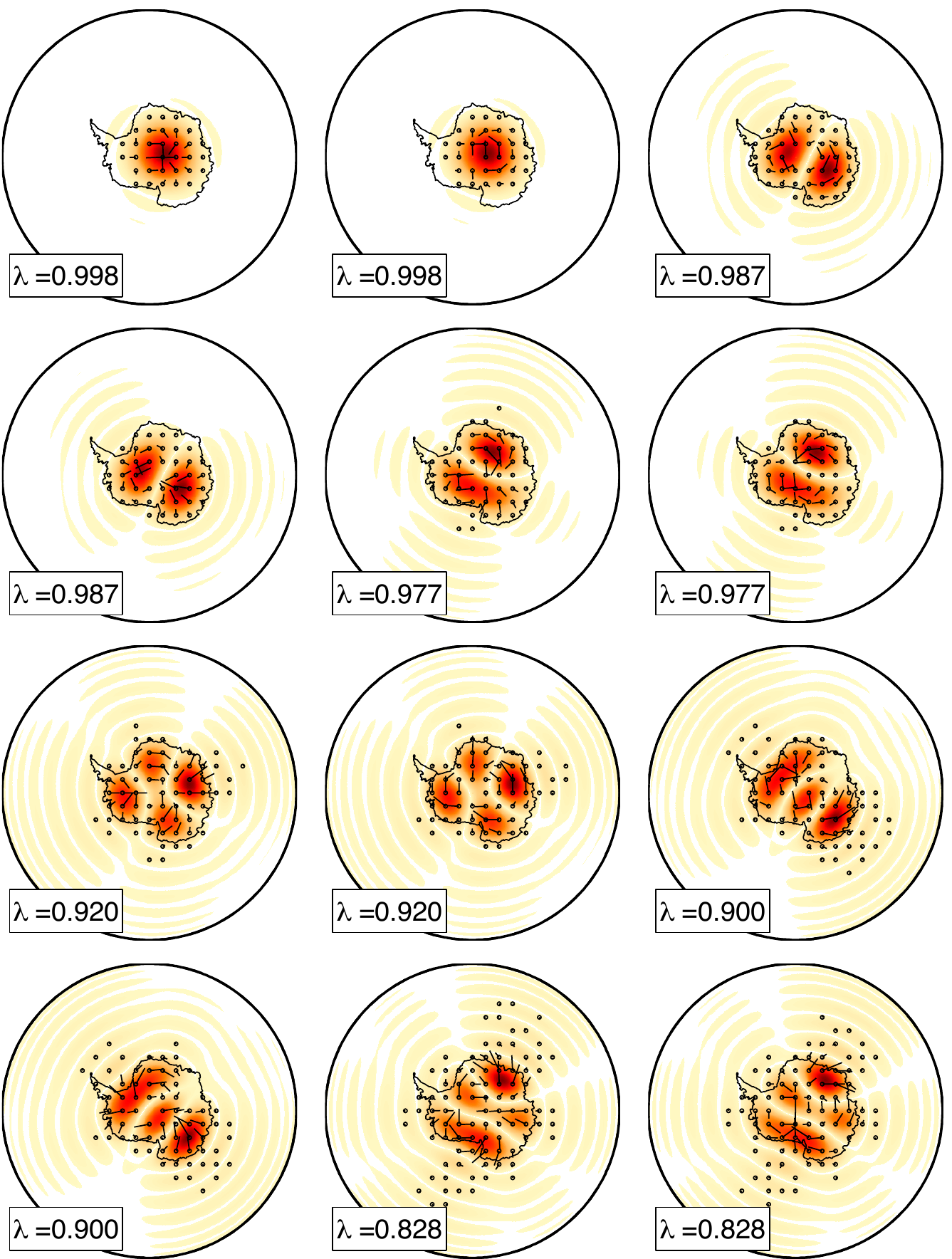}\figend
 \caption{\label{tangential antarctica} Bandlimited $L=18$ tangential
    eigenfields $\bg_1,\bg_2,\ldots,\bg_{12}$ that are optimally
    concentrated within Antarctica.  The concentration factors
    $\lambda_1,\lambda_2,\ldots,\lambda_{12}$ are indicated. The
    rounded tangential Shannon number is $N^t=20$; Format is identical
    to that in Figure~\ref{tangential africa}.}
\end{centering}
\end{figure}

Figures~\ref{tangential africa}--\ref{tangential antarctica} are map
views of the twelve best-concentrated tangential eigenfields
$\bg_1(\rvec),\bg_2(\rvec),\ldots,\bg_{12}(\rvec)$ for the continents
Africa and Antarctica, at $L=18$.  In either case, pairs of
eigenfields have identical absolute values and the same associated
eigenvalues but show vectorial directions that are pointwise
perpendicular, see Section~\ref{eigenvalue pairs} 
All eigenfields for Antarctica have a singularity or are zero at
the south pole. This comes from the fact that all tangential vector
spherical harmonics~(\ref{Blm})--(\ref{Clm}) either have a singularity
or are zero at the south pole. Both figures show that the first 12
tangential eigenfields are well concentrated, which is also reflected
by the tangential Shannon numbers, $N^t=42$ for Africa and $N^t=20$
for Antarctica. In both cases, the absolute values of the first two
eigenfields are roughly circular domes centered in the middle of each
continent.  Subsequent orthogonal eigenfields $\bg_3,\bg_4,\ldots$
exhibit lobes in previously uncovered regions. In
Figure~\ref{tangential africa}, West Africa begins to be reasonably
well covered by $\bg_5$ and $\bg_6$, while Southern Africa is
uncovered until $\bg_7$ and $\bg_8$.  Later, increasingly more
oscillatory eigenfields cover smaller geographical features.  For
Antarctica, the third and fourth eigenfields $\bg_3$ and $\bg_4$ begin
to resolve the South America-facing (western) and the Australia-facing
(eastern) part of Antarctica, while the fifth and the sixth
eigenfields $\bg_5$ and $\bg_6$ resolve the Africa-facing (northern)
region and the region around the Transantarctic Mountains. Subsequent
eigenfields show more nodal lines and resolve smaller geographical
features.

\begin{figure}
\begin{centering}
 \includegraphics[angle=0, width=1\textwidth]{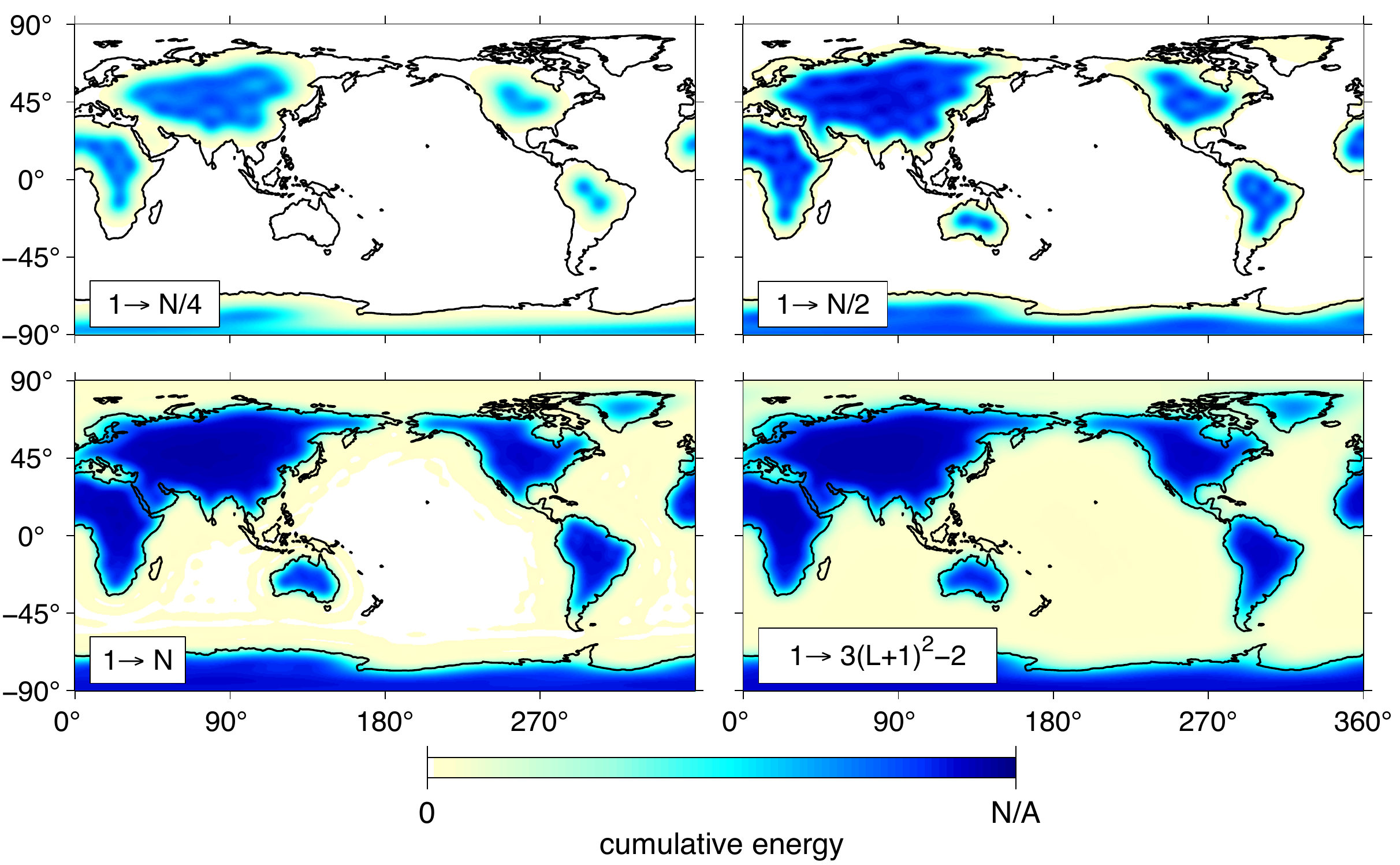}\figend
 \caption{\label{consum} Cumulative eigenvalue-weighted energy of the
   first $N/4,N/2,N$ and all $3(L+1)^2-2$ eigenfields that are
   optimally concentrated within the ensemble of Eurasia, Africa,
   North America, South America, Antarctica, Australia, and Greenland.
   The bandwidth is $L=18$; the cumulative fractional area is
   $A/(4\pi)=27.92$\%; the rounded Shannon number~$N=302$. The darkest
   blue on the color bar corresponds to the expected value
   (\ref{cumulative energy}) of the sum, as shown. Regions where the
   value is smaller than one hundredth of the $N/A$ are left white.}
\end{centering}
\end{figure}

Figure~\ref{consum} shows the eigenvalue-weighted sum of absolute
squares $\sum_\alpha \lambda_\alpha\lvert\bg_\alpha(\rvec)\rvert^2$ of
the $L=18$ bandlimited eigenfields of Earth's seven landmasses.  The
eigenfields $\bg_1(\rvec),\bg_2(\rvec),\ldots,\bg_{3(L+1)^2-2}(\rvec)$
can be found by diagonalizing the $[3(L+1)^2-2]\times[3(L+1)^2-2]$
matrix~(\ref{K matrix}) formed by summing the matrices
$\sK_\text{Eurasia}+\sK_\text{Africa}+ \ldots$ of each of the regions.
The fractional area covered by all seven regions combined is
$A/(4\pi)=27.92\%$, and the corresponding rounded Shannon
number~$N=302$; the figure shows the partial sums of the first
$N/4,N/2$ and $N$ terms, and the full sum of all $3(L+1)^2-2=1081$
terms.  It is apparent that the first~$N$ eigenfunctions uniformly
cover the target region; by adding the remaining, poorly concentrated,
$3(L+1)^2-2-N=779$ terms, we only marginally improve the coverage.
Because of its small size, Greenland does not appear until the $1
\rightarrow N/2$ partial sum. Even after the $1 \rightarrow N$ partial
sum, Greenland's coverage is not perfect, as  expected from its small
Shannon number ($N=5$ for $L=18$). 

\subsection{Spacelimited Fields}

As described in Section~\ref{spectral concentration}, the spatially
limited, spectrally concentrated vector fields~$\bh(\rvec)$ for a
region~$R$ and bandlimit~$L$ can be calculated by either spacelimiting
the spatially concentrated bandlimited fields for the same region~$R$
and the same bandlimit~$L$, as expressed by~(\ref{h normalization}),
or by multiplying the coefficient vector~$\sg$ of the spatially
concentrated band\-limited field with a rectangular kernel
matrix~$\sK$ of infinite bandwidth in the first dimension.

\begin{figure}
\begin{centering}
 \includegraphics[angle=0,
 width=0.875\textwidth]{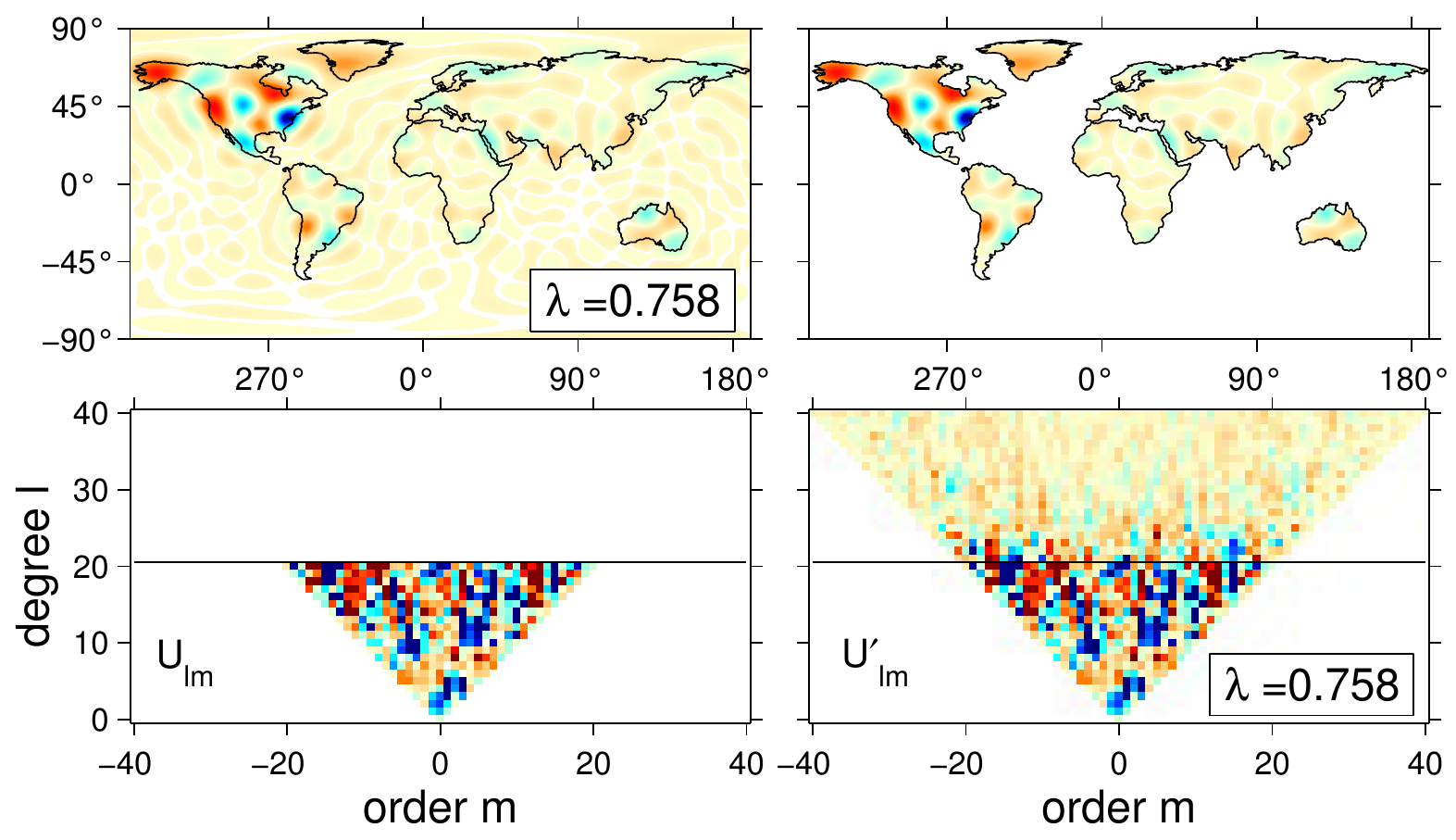}\figend 
 \caption{\label{spectral radial} The 80th best spatially concentrated
   bandlimited radial eigenfield ($\bg_{80}$, left panels) and the
   80th best spectrally concentrated spacelimited eigenfield
   ($\bh_{80}$, right panels) for a spatial domain which is the
   ensemble of Eurasia, Africa, North America, South America,
   Australia and Greenland, and a bandwidth $L=20$. The upper panels
   show the intensity and direction of the fields in the radial
   direction (blue inwards, red outwards). Regions where the absolute
   value is smaller than one hundredth of the maximum absolute value
   on the sphere are left white. The lower panels show the expansion 
   coefficients for the radial vector harmonics~$\Plm$ of the fields
   shown in the panels above (blue negative, red positive).}
\end{centering}
\end{figure}
\begin{figure}
\begin{centering}
 \includegraphics[angle=0,
 width=0.875\textwidth]{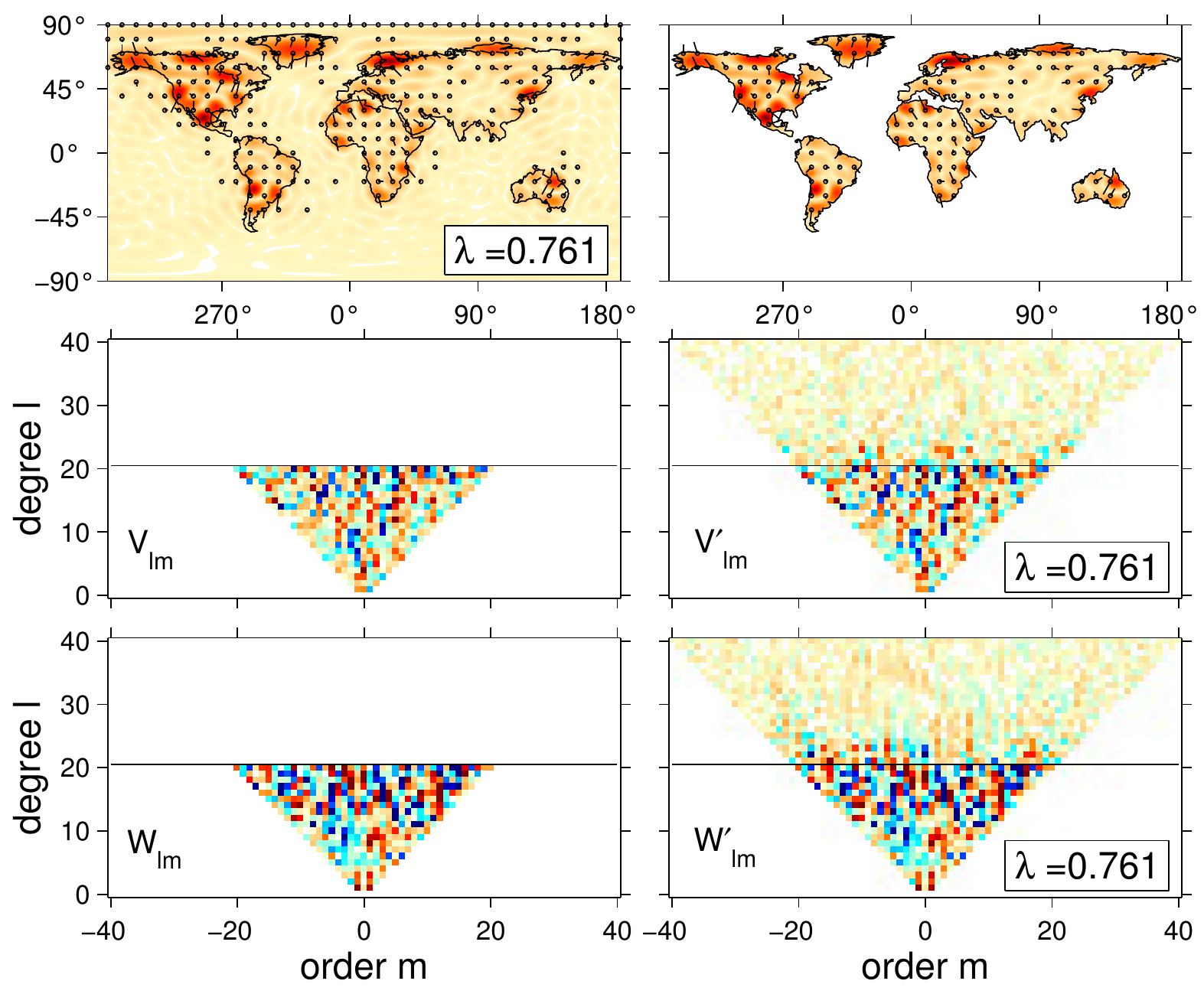}\figend 
 \caption{\label{spectral tangential} The 160th best spatially
   concentrated bandlimited tangential eigenfield~($\bg_{160}$, left
   panels) and the 160th best spectrally concentrated spacelimited
   eigenfield~($\bh_{160}$, right panels) for a spatial domain which
   is the ensemble of Eurasia, Africa, North America, South America,
   Australia and Greenland, and a bandwidth $L=20$.  Uppermost panels
   show the intensity and direction of the fields. Regions where the
   absolute value is smaller than one hundredth of the maximum
   absolute value on the sphere are left white.  Middle and lower
   panels show the expansion coefficients for the tangential vector
   harmonics~$\Blm$ and~$\Clm$, respectively, of the fields shown in
   the uppermost panels (blue negative, red positive).}
\end{centering}
\end{figure}

Figures~\ref{spectral radial} and~\ref{spectral tangential} show such
a construction for the combined six regions of Eurasia, Africa, North
and South America, Australia and Greenland and a bandlimit of~$L=20$.
Due to the block-diagonal shape of matrix $\sK$ in~(\ref{K matrix}),
the radial and tangential optimization problems are decoupled and were
solved independently. The upper left panel of Figure~\ref{spectral
  radial} shows the 80th best radial Slepian function $\bg(\rvec)$,
which by the measure~(\ref{radial coefficient optimization problem})
has $75.8\%$ of its energy within the target region. Blue stands for
inwards and red for outwards-pointing vectors. Areas with intensity of
less than one percent of the maximum value are left white. The lower
left panel of this figure shows the spherical-harmonic coefficients
$\sg$ of this radial Slepian field.  Shades of blue denote negative
coefficient values and red positive values. Due to the bandlimitation,
all coefficients with degree higher than $L=20$ are zero. The upper
right panel of Figure~\ref{spectral radial} shows the spatially
truncated radial Slepian field $\bh(\rvec)$ and the lower right panel
its spherical-harmonic coefficients~$\sh$. The coefficients~$\Ulm$ are
only shown up to $l=40$ but are nonzero to $l=\infty$ since
$\bh(\rvec)$ is perfectly spacelimited. The ratio~(\ref{spectral
  concentration equation}), of the energy in the coefficients below
$l=20$ to the total energy, is once again~$75.8\%$, illustrating the
equivalence between the spatial and spectral concentration problems.

Figure~\ref{spectral tangential} illustrates the same procedure
applied to tangential fields. The upper left panel shows the 160th
best spatially concentrated bandlimited tangential field, with maximal
spherical-harmonic degree $L=20$. The middle left panel and the bottom
left panel show the vector spherical-harmonic coefficients $\Vlm$ and
$\Wlm$, respectively. Again, the vector spherical-harmonic
coefficients at degrees above the bandlimit~$L=20$ are zero. The right
panels show the spacelimited and spectrally concentrated tangential
vector Slepian field constructed from the bandlimited spatially
concentrated vector Slepian field shown on the left, both in their
spatial (uppermost panel) and spectral (middle and lower panels)
renditions. 

\subsection{Constructive Approximation}

\begin{figure}
\begin{centering}
\includegraphics[angle=0,
width=0.75\textwidth]{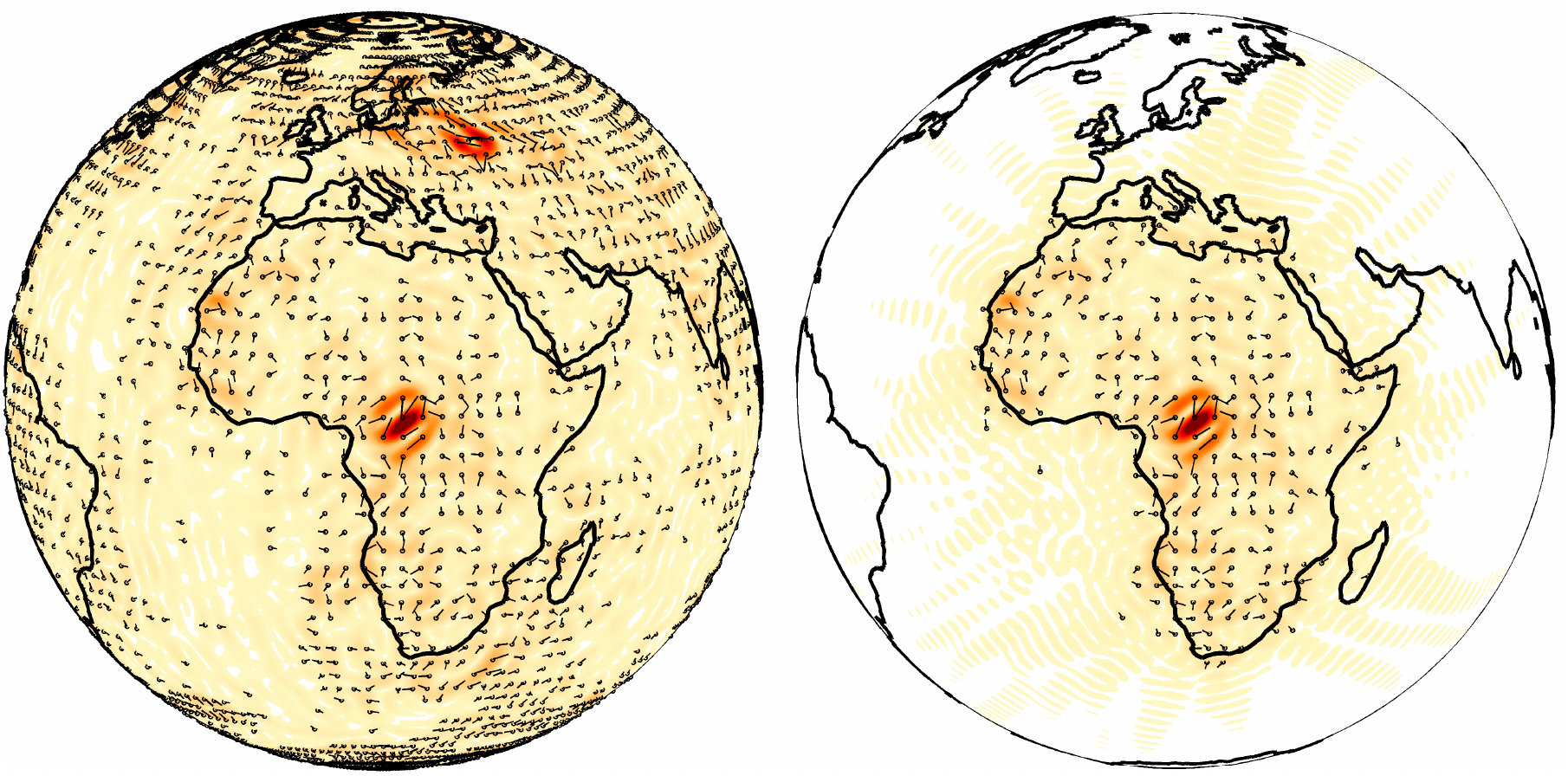}\\[1em] 
\includegraphics[angle=0,
width=0.635\textwidth]{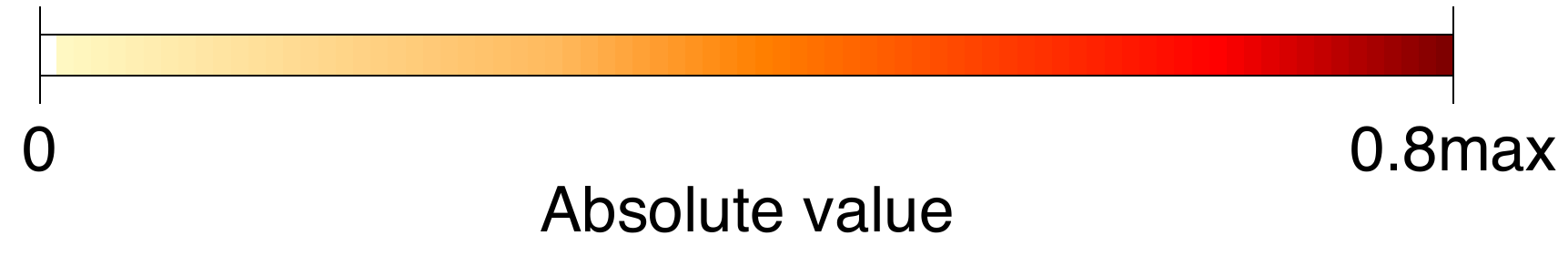}\figend
\caption{\label{magfield} A tangential geophysical vector field (left)
  and its reconstruction (right) using vectorial Slepian functions
  designed to maximize their spatial concentration over Africa. The
  bandlimit for both the original field and the Slepian basis is $L =
  72$. There are 10,656 vectorial basis functions in the original
  field, and the same number of Slepian functions from which to choose
  for the reconstruction. The Shannon number $N^t=620$. The bottom
  panel shows a reconstruction using the $924=1.5 N^t$
  best-concentrated Slepian functions for Africa. The error and bias
  over Africa, as defined in~(\ref{error and bias}), are 0.4\% and
  14\%, respectively.}
\end{centering}
\end{figure}

Finally, in order to demonstrate the spatial focusing capabilities of
the bandlimited, spatially concentrated vector Slepian fields for an
actual data example, we reconstruct a global tangential vector field,
$\bu$, by approximating it with fields $\bv_J$ that use an increasing
number, $J$, of tangential vector Slepian functions:
\begin{equation}\label{bvj}
  \bv_J = \sumalphaJ u_\alpha\hsom \bga,
\end{equation} 
The coefficients $u_\alpha$ are obtained by forming the inner product
of the input field $\bu$ with the $\alpha$ best-concentrated vector
Slepian functions $\bga$. We define the relative error~$\epsilon_J$
over the domain, and the leakage~$b_J$ to its complement, by 
\begin{equation}\label{error and bias}
\epsilon_J=\sqrt\fracd{\lVert\bu-\bv_J\rVert_R^2}{\lVert\bu\rVert_R^2}
\also
b_J=\sqrt\fracd{\lVert\bv_J\rVert^2_{\Omega\setminus
    R}}{\lVert\bu\rVert^2_{\Omega\setminus R}}, 
\end{equation}
which we will use to assess the performance of the reconstruction.
For bandlimited tangential fields~$\bu$ the measure~$\lVert\bu
\rVert_R^2$ can be calculated using the matrix~$\sQ$ from~(\ref{Q
  matrix}) and the vector of expansion coefficients
$\su=(\ldots,\Ulm,\ldots,\Vlm,\ldots,\Wlm,\ldots)^\iT$ by evaluating
the expression~$\lVert\bu \rVert_R^2=\su^\sT \sQ\hsom\su$.  The error
decreases with increasing number of Slepian functions~$J$. The bias
increases with~$J$. Our goal is to obtain a small reconstruction error
within the region~$R$ while simultaneously keeping the outside leakage
bias small.

Figure~\ref{magfield} shows the outcome of such an experiment
conducted on the terrestrial crustal-field model NGDC-720~V3
\cite{Maus2010}. We multiply the spherical-harmonic coefficients with
the corresponding~$\Blm$ vector harmonics up to bandlimit~$L=72$. The
left panel of Figure~\ref{magfield} shows the tangential vector field
that results: this is used as our input. The right panel shows the
reconstruction using the $1.5N=924$ best-concentrated tangential
vector-Slepian fields for Africa and the same bandlimit~$L=72$.
  While we chose $1.5N$ here for convenience, in real-world
  applications the optimal choice for the Slepian truncation would
  depend on the behavior of the signal-to-noise ratio of the
  data~\cite{Simons+2006b,Schachtschneider+2010,Schachtschneider+2012}.
The relative error and bias of the reconstruction
over Africa, as defined by~(\ref{error and bias}), are~0.4\% and~14\%,
respectively.

\begin{figure}
\begin{centering}
 \includegraphics[angle=0, width=0.6\textwidth]{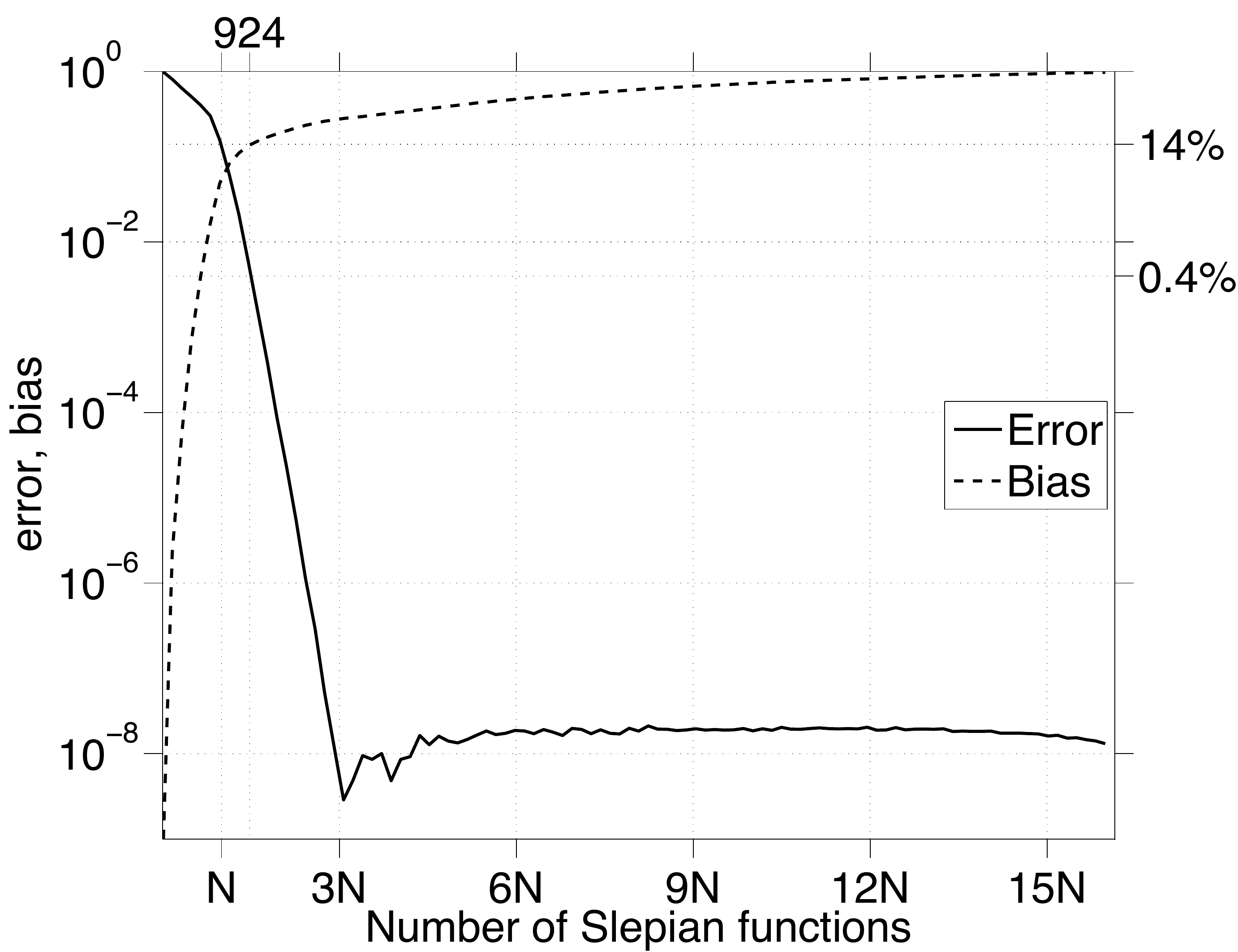}\figend
 \caption{\label{magfield error} Reconstruction error and bias over
   Africa as defined in~(\ref{error and bias}), versus the number of
   vector Slepian functions used to describe the global vector field
   as shown in Figure~\ref{magfield}, quoted as a multiple of the
   Shannon number for this problem, $N^t=620$.}
\end{centering}
\end{figure}

Figure~\ref{magfield error} shows the evolution of error and bias for
reconstructions using different numbers of Slepian-field terms in the
expansion~(\ref{bvj}). The more Slepian fields are being used, the
smaller the error over Africa, but the larger the leakage into the
complimentary region outside of Africa. The relative reconstruction
error over Africa drops quickly and reaches numerical noise level
after $J=1800$ Slepian function terms. 

\section{Conclusion}

It is possible to construct for the unit sphere a regionally optimally
concentrated orthogonal family of bandlimited vector
spherical-harmonic fields by solving either a Fredholm integral
eigenvalue problem in the spatial domain, or, equivalently in the
spectral domain solving a symmetric finite-dimensional matrix
eigenvalue problem. The eigenvalues $0<\lambda<1$ are measures of the
spatial concentrations of their corresponding bandlimited vector
fields $\bg(\rvec)$ and spectral concentration of the spacelimited
eigenfields $\bh(\rvec)$, constructed from the $\bg(\rvec)$ by setting
the values to zero outside of the target region.  The full vectorial
problem decomposes into independent radial and tangential parts. The
radial problem is equivalent to the scalar spherical spatiospectral
optimization problem~\cite{Simons+2006a}.  The number of
well-concentrated radial eigenfields is $N^r=(L+1)^2A/(4\pi)$ and the
number of well-concentrated tangential eigenfields is
$N^t=[2(L+1)^2-2]A/(4\pi)$. Here~$L$ denotes the bandwidth and $A$ the
area of the target region. The Shannon numbers $N^r$ and $N^t$ can be
interpreted as the dimensions of the spaces of radial vector fields
$\bg^r(\rvec)$, or tangential vector fields $\bg^t(\rvec)$
respectively, that can be simultaneously concentrated within a
subregion~$R$ of the sphere and within a spectral interval $0\leq
l\leq L$. In the special case of a circular polar cap, the kernel
matrices can be computed analytically and decomposed into smaller
eigenvalue problems.

%
%

Vectorial Slepian functions on the sphere are an emerging tool for the
analysis and representation of essentially space- and bandlimited
vector-valued functions on the surface of the unit sphere. In this
contribution we have described their construction, shown various
examples, and suggested their use in the constructive approximation of
vectorial signals on the sphere, as may arise, for instance, in the
fields of geophysics, planetary science, medical imaging and optics,
where prior work has previously considered a number of special cases
of the vectorial concentration on the sphere~\cite{Mitra2006,Jahn2012}
that we have treated more completely here.

  The ability of scalar Slepian functions on the sphere to perform
  localized bandlimited analysis has led to observations made from
  global data that remain obscured when applying global spherical
  harmonic analysis. For example, changes in local gravity after the
  2004 Sumatra earthquake were detected~\cite{Han2008d,Simons+2009b}
  and shown to be invisible via a global spherical harmonic analysis
  of the same data. Similarly, in analyzing global gravity data, the
  potential of scalar Slepian functions to detect local ice mass
  changes over Greenland was clearly demonstrated~\cite{Harig+2012}.
  Judging from the equivalence in properties between the vectorial
  Slepian functions and the scalar Slepian functions in multiple
  Cartesian and spherical dimensions,
it is likely that the impact of vectorial spherical Slepian functions
on multidimensional vectorial signal processing will be as profound as
the classical prolate spheroidal wave functions have been, and
continue to be, in the study of time series, and this in a wide
variety of scientific and engineering fields.

\section*{Acknowledgments}

A.~P. and F.J.S~thank Volker Michel for valuable discussions, and the
Ulrich Schmucker Memorial Trust and the Swiss National Science
Foundation for financial support.  F.J.S.~thanks Tony Dahlen
(1942--2007) and Liying Wei for their contributions to the material
developed in this manuscript. We thank the anonymous referee and the
Editor-in-Chief Charles Chui. This work was partially supported by
National Science Foundation grants EAR-1014606 and EAR-1150145 to
F.J.S. To the historians of science we like to mention that we
discovered reference~\cite{Jahn2012} only after concluding our
analysis, when completing our manuscript for submission.

\bibliographystyle{ieeetr}
\bibliography{slepian,/u/fjsimons/BIBLIO/bib}

\end{document}